\pdfoutput=1
\documentclass[tightenlines, twocolumn, notitlepage, nofootinbib, superscriptaddress]{revtex4-2}

\usepackage[dvipsnames]{xcolor}
\usepackage{amsmath}
\usepackage{amssymb}
\usepackage{bm}
\usepackage{bibunits}
\usepackage{changes}
\usepackage{enumerate}
\usepackage{floatrow}
\usepackage{hyperref}
\usepackage[symbol]{footmisc}
\usepackage{graphicx}
\usepackage{units}
\usepackage{xspace}
\usepackage{booktabs}
\usepackage{listings}
\usepackage[frozencache=true,cachedir=minted-cache]{minted} 
\usepackage{multirow}
\usepackage{gensymb}
\usepackage{algorithm2e}
\usepackage[utf8x]{inputenc}
\usepackage[strings]{underscore}
\usepackage[english]{babel}
\RestyleAlgo{ruled}

\usepackage{soul}

\newfloatcommand{capbtabbox}{table}[][\FBwidth]
\usepackage{blindtext}

\graphicspath{{./}, {figures/}}
\setcitestyle{super}

\newcommand{\olympus}{\textsc{Olympus}\xspace}
\newcommand{\gpyopt}{\textsc{GPyOpt}\xspace}
\newcommand{\deap}{\textsc{DEAP}\xspace}
\newcommand{\phoenics}{\textsc{Phoenics}\xspace}
\newcommand{\gryffin}{\textsc{Gryffin}\xspace}
\newcommand{\dragonfly}{\textsc{Dragonfly}\xspace}
\newcommand{\chimera}{\textsc{Chimera}\xspace}

\DeclareMathOperator*{\argmin}{arg\,min}

\makeatletter
\renewcommand*{\p@subsection}{\thesection.}
\makeatother

\usepackage{booktabs, arydshln}

\SetKwComment{Comment}{/* }{ */}

\makeatletter
\renewcommand*{\p@subsection}{\thesection.}
\makeatother

\begin{document}


	\title{\large{Bayesian optimization with known experimental and design constraints for chemistry applications}}

	\date{\today}
	
	\author{Riley J. Hickman}
	\thanks{These authors contributed equally}
	\affiliation{Chemical Physics Theory Group, Department of Chemistry, University of Toronto, Toronto, ON M5S 3H6, Canada}
	\affiliation{Department of Computer Science, University of Toronto, Toronto, ON M5S 3G4, Canada}
	
	\author{Matteo Aldeghi}
	\thanks{These authors contributed equally}
	\affiliation{Chemical Physics Theory Group, Department of Chemistry,
	University of Toronto, Toronto, ON M5S 3H6, Canada}
	\affiliation{Department of Computer Science, University of Toronto, Toronto, ON M5S 3G4, Canada}
	\affiliation{Vector Institute for Artificial Intelligence, Toronto, ON M5S 1M1, Canada}
	\affiliation{Department of Chemical Engineering, Massachusetts Institute of Technology, Cambridge, MA 02139, United States}
	
	\author{Florian H\"ase}
	\affiliation{Chemical Physics Theory Group, Department of Chemistry, University of Toronto, Toronto, ON M5S 3H6, Canada}
	\affiliation{Department of Computer Science, University of Toronto, Toronto, ON M5S 3G4, Canada}
	\affiliation{Vector Institute for Artificial Intelligence, Toronto, ON M5S 1M1, Canada}
	\affiliation{Department of Chemistry and Chemical Biology, Harvard University, Cambridge, MA 02138, United States}

	\author{Al\'an Aspuru-Guzik}
	\email{alan@aspuru.com}
	\affiliation{Chemical Physics Theory Group, Department of Chemistry, University of Toronto, Toronto, ON M5S 3H6, Canada}
	\affiliation{Department of Computer Science, University of Toronto, Toronto, ON M5S 3G4, Canada}
	\affiliation{Vector Institute for Artificial Intelligence, Toronto, ON M5S 1M1, Canada}
	\affiliation{Department of Chemical Engineering \& Applied Chemistry, University of Toronto, Toronto, ON M5S 3E5, Canada}
	\affiliation{Department of Materials Science \& Engineering, University of Toronto, Toronto, ON M5S 3E4, Canada}
	\affiliation{Lebovic Fellow, Canadian Institute for Advanced Research, Toronto, ON M5G 1Z8, Canada}


\begin{abstract}
Optimization strategies driven by machine learning, such as Bayesian optimization, are being explored across experimental sciences as an efficient alternative to traditional design of experiment. When combined with automated laboratory hardware and high-performance computing, these strategies enable next-generation platforms for autonomous experimentation. 
However, the practical application of these approaches is hampered by a lack of flexible software and algorithms tailored to the unique requirements of chemical research. One such aspect is the pervasive presence of constraints in the experimental conditions when optimizing chemical processes or protocols, and in the chemical space that is accessible when designing functional molecules or materials. Although many of these constraints are known \textit{a priori}, they can be interdependent, non-linear, and result in non-compact optimization domains. 
In this work, we extend our experiment planning algorithms \textsc{Phoenics} and \textsc{Gryffin} such that they can handle arbitrary known constraints via an intuitive and flexible interface. We benchmark these extended algorithms on continuous and discrete test functions with a diverse set of constraints, demonstrating their flexibility and robustness. In addition, we illustrate their practical utility in two simulated chemical research scenarios: the optimization of the synthesis of \textit{o}-xylenyl Buckminsterfullerene adducts under constrained flow conditions, and the design of redox active molecules for flow batteries under synthetic accessibility constraints.
The tools developed constitute a simple, yet versatile strategy to enable model-based optimization with known experimental constraints, contributing to its applicability as a core component of autonomous platforms for scientific discovery.

\end{abstract}

	\maketitle
	
	\begin{bibunit}[ieeetr]


\section{Introduction} \label{sec:introduction}

The design of advanced materials and functional molecules often relies on combinatorial, high-throughput screening strategies enabled by high-performance computing and automated laboratory equipment. Despite the successes of high-throughput experimentation in chemistry,~\cite{mcnally_discovery_2011,collins_contemporary_2014} biology,~\cite{blay_high-throughput_2020,zeng_high-throughput_2020} and materials science,~\cite{cheng_accelerating_2015} these approaches typically employ exhaustive searches that scale exponentially with the size of the search space. Data-driven strategies that can adaptively search parameter spaces without the need for exhaustive exploration are thus replacing traditional design of experiment approaches in many instances. These strategies use machine-learnt surrogate models trained on all data generated through the experimental campaign, and are updated each time new data is collected. One such approach is Bayesian optimization which, based on the surrogate model, defines a utility function that prioritize experiments based on their expected informativeness and performance.~\cite{Mockus:1975,Mockus:1978,Mockus:2012} These data-driven optimization strategies have already demonstrated superior performance in chemistry and materials science applications, e.g., in reaction optimization,\cite{shields_bayesian_2021,christensen_data-science_2021} the discovery of magnetic resonance imaging agents,\cite{Reis:2021} the fabrication of organic photovoltaic materials,\cite{Langner:2020,MacLeod:2020} virtual screening of ultra-large chemical libraries,\cite{Graff:2021} and the design of mechanical structures with additively manufactured components.\cite{Gongora:2020}

Machine learning-driven experiment planning strategies can also be combined with automated laboratory hardware or high-performance computing to create self-driving platforms capable of achieving research goals autonomously.~\cite{hase_next-generation_2019,roch_chemos_2020,correa-baena_accelerating_2018,sstein_progress_2019,stach_autonomous_2021,coley_autonomous_2020,coley_autonomous_2020-1} Prototypes of these autonomous research platforms have already shown promise in diverse applications, including the optimization of chemical reaction conditions,~\cite{shields_bayesian_2021,christensen_data-science_2021} the design of photocatalysts for the production of hydrogen from water,~\cite{burger_mobile_2020} the discovery of battery electrolytes,~\cite{dave_autonomous_2020} the design of nanoporous materials with tailored adsorption properties,~\cite{deshwal_bayesian_2021} the optimization of multicomponent polymer blend formulations for organic photovoltaics,~\cite{Langner:2020} the discovery of phase-change memory materials for photonic switching devices,~\cite{kusne_--fly_2020} and self-optimization of metal nanoparticle synthesis,~\cite{tao_self-driving_2021} to name a few.~\cite{wigley_fast_2016,MacLeod:2020} While self-driving platforms seem poised to deliver a next-generation approach to scientific discovery, their practical application is hampered by a lack of flexible software and algorithms tailored to the unique requirements of chemical research.

To provide chemistry-tailored data-driven optimization tools, our group has developed \phoenics\cite{hase_phoenics_2018} and \gryffin,\cite{hase_gryffin_2021} among others\cite{hase_chimera_2018,Aldeghi:2021_golem,Hickman:2021_gemini}. \phoenics is a linear-scaling Bayesian optimizer for continuous spaces that uses a kernel regression surrogate model and natively supports batched optimization. \gryffin is an extension of this algorithm to categorical, as well as mixed continuous-categorical spaces. Furthermore, \gryffin is able to leverage expert knowledge in the form of descriptors to enhance its optimization performance, which was found particularly useful in combinatorial optimizations of molecules and materials.\cite{hase_gryffin_2021} As \gryffin is the more general algorithm, and \phoenics is included within its capabilities, from here on we will refer only to \gryffin. These algorithms have already found applications ranging from the optimization of reaction conditions\cite{christensen_data-science_2021} and synthetic protocols,\cite{seifrid_routescore_2021,tao_self-driving_2021} to that of manufacturing of thin film materials\cite{MacLeod:2020} and organic photovoltaic devices.\cite{Langner:2020} However, a number of extensions are still required to make these tools suitable to the broadest range of chemistry applications. In particular, \gryffin, like the majority of Bayesian optimization tools available, does not handle known experimental or design constraints.

There are often many constraints on the experiment being performed or molecule being designed. A flexible data-driven optimization tool should be able to accommodate and handle such constraints. The type of constraints typically encountered may be separated into those that affect the objectives of the optimization (e.g., reaction yield, desired molecular properties), and those that affect the optimization parameters (e.g., reaction conditions). Those affecting the objectives usually arise in multi-objective optimization, where one would like to optimize a certain property while constraining another to be above/below a desired value.\cite{walker_tuning_2017,Gelbart:2014,hase_chimera_2018} For instance, we might want to improve the solubility of a drug candidate, while keeping its protein-binding affinity in the nanomolar range. Conversely, parameter constraints limit the range of experiments or molecules we have access to. Depending on the source of the constraints, these may be referred to as \textit{known} or \textit{unknown}. Known constraints are those we are aware of \textit{a priori},~\cite{gramacy_optimization_2010,gelbart_bayesian_2014,ariafar_admmbo_2019,antonio_sequential_2021} while unknown ones are discovered through experimentation~\cite{gramacy_optimization_2010,gelbart_bayesian_2014,ariafar_admmbo_2019,antonio_sequential_2021}. For instance, a known constraint might enforce the total volume for two liquids to not exceed the available volume in the container in which they are mixed. While this poses a restriction on the parameter space, we are aware of it in advance and can easily compute which regions of parameter space are infeasible. An unknown constraint may instead be the synthetic accessibility in a molecular optimization campaign. In this case, we might not know in advance which areas of chemical space are easily accessible, and have to resort to trial and error to identify feasible and infeasible synthetic targets. While constraints of the objectives were the subject of our previous work\cite{hase_chimera_2018}, and unknown constraints of the parameters are the subject of on-going work, this paper focuses on data-driven optimization with known parameter constraints, which we will refer to simply as \textit{known constraints} from here on.

Generally, known constraints arise due to physical or hardware restrictions, safety concerns, or user preference. An example of a physically imposed constraint is the fact that the temperature of a reaction cannot exceed the boiling temperature of its solvent. As such, one may want temperature to be varied in the interval $10 < T < 100$\textdegree C  for experiments using water, and in $10 < T < 66$\textdegree C for experiments using tetrahydrofuran. The fact that the sum of volumes of different solutions cannot exceed that of the container they are mixed in is an example of a hardware-imposed constraint. In synthetic chemistry, specific combinations of reagents and conditions might need to be avoided for safety reasons instead. Finally, constraints could also be placed by the researchers to reflect prior knowledge about the performance of a certain protocol. For example, a researcher might know in advance that specific combinations of solvent, substrate, and catalyst will return poor yields. These examples are not natively handled by \gryffin and the majority of data-driven optimization tools currently available. In fact, given any number of continuous or categorical parameters, their full Cartesian product is assumed to be accessible by the optimization algorithm. Returning to the example where solvents have different boiling temperatures, this means that if the optimization range for the variable $T$ is set to $10-100$\textdegree C, this range will be applied to all solvents considered. In practice, known constraints are often interdependent, non-linear, and can result in non-compact optimization domains. 

In this work, we extend the capabilities of \gryffin to optimization over parameter domains with known constraints. First, we provide a formal introduction to the known constraint problem and detail how \gryffin was extended to be able to flexibly handle such scenarios. Then, we benchmark our new constrained version of \gryffin on a range of analytical functions subject to a diverse set of constraints. Finally, we demonstrate our method on two chemistry applications: the optimization of the synthesis of \textit{o}-xylenyl Buckminsterfullerene adducts under constrained flow conditions, and the design of redox active molecules for flow batteries under synthetic accessibility constraints. Across all tests, we compare \gryffin's performance to that of other optimization strategies, such as random search and genetic algorithms, which can also handle complex constraints.


\section{Methods} \label{sec:methods}

An optimization task involves the identification of parameters, $\bm{x}$, that yield the most desirable outcome for an objective $f(\bm{x})$. In a chemistry context, these parameters may be experimental conditions or different R groups in a molecule, while the objectives may be the yield of a reaction or absorbance at a specific wavelength. Formally, for a minimization problem, the solution of the optimization is the set of parameters that minimizes the objective $f(\bm{x})$,

\begin{align*}
  \bm{x}^*  = \argmin_{\bm{x} \in \mathcal{X}} f (\bm{x}) \text{,}
\end{align*}
    
where $\mathcal{X}$ is the optimization domain, or parameter space; i.e., the space of all experimental conditions that could have been explored during the optimization. In a Bayesian optimization setting, the objective function $f$ is considered to be unknown, but can be empirically evaluated at specific values of $\bm{x}$. Evaluating $f(\bm{x})$ is assumed to be expensive and/or time consuming, and its measurement subject to noise. We also assume that we have no access to gradient information about $f$. 

In a \textit{constrained} optimization problem, $f$ can be evaluated only for a subset of the optimization domain, $\mathcal{C} \subset \mathcal{X}$. A constraint function $c(\bm{x})$ determines which parameters $\bm{x}$ are feasible and thus in $\mathcal{C}$, and which are not. Contrary to \textit{unknown} constraints, known constraints are those for which we have access to $c(\bm{x})$ \textit{a priori}, that is, we are aware of them before performing the optimization. Importantly, we can evaluate the functions $c(\bm{x})$ and $f(\bm{x})$ independently of one another. The solution of this constrained optimization problem may be written formally as

\begin{align*}
  \bm{x}^*  = \argmin_{\bm{x} \in \mathcal{X}} f (\bm{x}) \text{,}\\
  \text{s.t.}\ c(\bm{x}) \mapsto \text{feasible} \text{.}
\end{align*}

The physically meaning of $c(\bm{x})$ is case dependent and domain specific. In chemistry, known constraints may reflect safety concerns, physical limits imposed by laboratory equipment, or simply researcher preference. The types of known constraints described above, and the subject of this work, are \textit{hard} constraints. They restrict $\mathcal{X}$ irrefutably and no measurement outside of $\mathcal{C}$ is permitted. However, the use of \textit{soft} constraints has also been explored. While soft constraints bias the optimization algorithm away from regions that are thought to yield undesirable outcomes, they ultimately still allow the full exploration of $\mathcal{X}$. As such, soft constraints have been used as means to introduce inductive biases into an optimization campaign based on prior knowledge.~\cite{sun_physical_2020,liu_machine_2021}

The goal of this work is to equip \gryffin with the ability to satisfy arbitrary constraints, as described by a user-defined constraint function $c(\bm{x})$. This can be achieved by constraining the optimization of the acquisition function $\alpha$.\cite{Shahriari:2016} $\alpha(\bm{x})$ defines the utility (or risk) of evaluating the parameters $\bm{x}$, and the parameters proposed by the algorithm are those that optimize $\alpha(\bm{x})$. Thus, constrained optimization of the acquisition function also constrains the overall optimization problem accordingly. However, contrary to objective function $f$, the acquisition function $\alpha$ is easy to optimize, as its analytical form is known and can be evaluated cheaply.

\begin{figure}[htb]
    \centering
    \includegraphics[width=0.95\columnwidth]{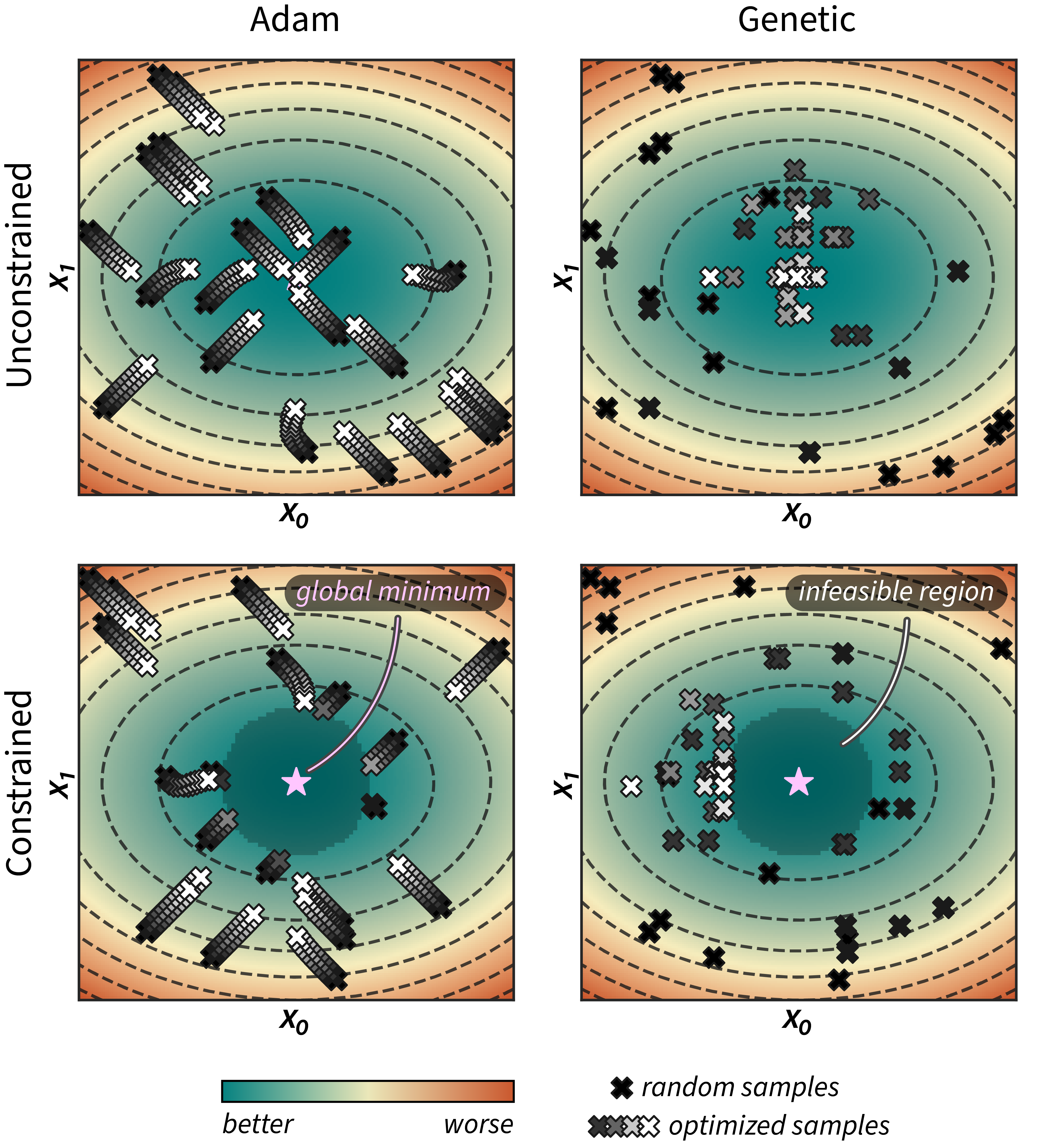}
    \label{fig:toy_illustration}
    \caption{Illustration of unconstrained (top row) and constrained (bottom row) acquisition function optimization using strategies based on the Adam optimizer and a genetic algorithm. Initial, random samples are shown as black crosses. Grey crosses represent the updated parameter locations for these initial samples, while white crosses show the final locations after ten optimization iterations. The purple star indicates the global optima of the unconstrained acquisition function, which lies in the infeasible region in the constrained example.}
    \label{fig:acq_demonstration}
\end{figure}

\subsection{Acquisition optimization in \textsc{Gyrffin}} \label{subsec:acq_optim_gryffin}

\gryffin's acquisition function, which is to be minimized, is defined as

\begin{align} \label{eq:gryffin_acq}
  \alpha(\bm{x}) = \frac{\sum_{k=1}^{n} f_k p_k(\bm{x}) + \lambda p_{\text{uniform}}(\bm{x})}{\sum_{k=1}^{n} p_k(\bm{x}) + p_{\text{uniform}}(\bm{x})} \text{,}
\end{align}

where $p_k(\bm{x})$ are the kernels of the kernel regression model used as the surrogate, $f_k$ are the measured objective function values, and $\lambda$ is a user-defined parameter that determines the exploration-exploitation behavior of the algorithm. The index $k$ refers to each past observation, for a total number of observations $n$.

For continuous parameters, \gryffin uses Gaussian kernels with prior precision $\tau$ sampled from a Gamma distribution, $\tau \sim \Gamma (a,b)$, with prior hyperparameters $a = 12n^2$ and $b = 1$. For categorical and discrete parameters, it uses Gumbel-Softmax\cite{Jang:2017,Maddison:2017} kernels with a prior temperature parameter of $0.5 + 10n^{-1}$. These prior parameters are, however, updated by a Bayesian neural network in light of the observed data.\cite{hase_phoenics_2018} As the prior precision of the kernels increases with the number of observations, the surrogate model is encouraged to fit the observations more accurately as more data is acquired.

Acquisition function optimization in \gryffin generally follows a two-step strategy in which a global search is followed by local refinement. First, sets of input parameters $\bm{x}_i$ are sampled uniformly from the optimization domain. By default, the number of samples is set to be directly proportional to the dimensionality of the optimization domain. Then, continuous parameters are optimized with a gradient method for a pre-defined number of iterations. While early versions of \gryffin employed second order methods such as L-BFGS, the default gradient-based approach is now Adam.~\cite{kingma_adam_2017} Categorical and discrete parameters are optimized instead following a hill-climbing approach (which we refer to as \textit{Hill}), in which each initial sample $\bm{x}_i$ is modified uniformly at random, one dimension at a time, and each move is accepted only if it improves $\alpha(\bm{x}_i)$.

In addition to gradient-based optimization of the acquisition function, in this work we have implemented gradient-free optimization via a genetic algorithm. The population of parameters is firstly selected uniformly at random, as described above. Then, each $\bm{x}_i$ in the population is evolved via crossover, mutation and selection operations, with elitism applied. This approach handles continuous, discrete, and categorical parameters by applying different mutations and cross-over operations depending on the type. A more detailed explanation of the genetic algorithm used for acquisition optimization is provided in SI Sec.~\ref{sisubsubsec:genetic_acquisition}. This approach is implemented in Python and makes use of the \deap library.~\cite{Fortin:2012,DeRainville:2012}

Figure~\ref{fig:acq_demonstration} provides a visual example of how these two approaches optimize the acquisition function of \gryffin. While Adam optimizes each random sample $\bm{x}_i$ (black crosses) via small steps towards better $\alpha(\bm{x}_i)$ values (grey to white crosses), the genetic approach does so via larger stochastic jumps in parameter space.

\subsection{Constrained acquisition optimization with Adam or Hill} \label{subsec:const_adam}

To constrain the optimization of the acquisition function according to a user-defined constraint function $c(\bm{x})$, we first sample a set of feasible parameters $\bm{x}_i$ with rejection sampling. That is, we sample $\bm{x} \sim \mathcal{X}$ uniformly from the optimization domain and we retain only samples that satisfy the constraint function $c(\bm{x})$. Sampling is performed until the desired number of feasible samples is drawn. Local optimization of parameters is then performed with Adam, as described above, but it is terminated as soon as an update results in $c({\bm{x}_i}) \mapsto \text{infeasible}$ (Fig.~\ref{fig:acq_demonstration}). For categorical and discrete parameters, Hill is used rather than Adam, but the constraint protocol is equivalent. In this case, after rejection sampling, any move in parameter space is accepted only if it improves $\alpha(\bm{x}_i)$ and subject to $\bm{x}_i \in \mathcal{C}$.

In addition to the above, we also modify the prior precisions of the kernels used by the surrogate model. In particular, we substitute the number of observations $n$ with the observation density $\rho = n / V_{\mathcal{C}}$ in the prior hyperparameters for the Gaussian and Gumbel-Softmax kernels, where $V_{\mathcal{C}}$ is the volume of the feasible region as a fraction of the overall optimization domain. When no constraints are used, $V_{\mathcal{C}}=1$, and the formulation reduces to the original one used by \gryffin. However, when known constraints are present, this approach ensures that the kernels' bandwidths reflect the information density of the data and avoids underfitting. As $V_{\mathcal{C}}$ is known and can be computed in advance, the user is asked to  provide it when using constraints. However, for convenience, and in case of fairly complex sets of constraints, we extend \gryffin to also numerically estimate $V_{\mathcal{C}}$. While this approach to kernel scaling may not be optimal in the presence of non-isotropic feasible regions, in which the viable optimization domain may expand to different extent in different dimensions, it is simple as well as independent of the details of the constraints, such that it results in no additional overheads.

\subsection{Constrained acquisition optimization with a genetic algorithm} \label{subsec:const_genetic}

The population of the genetic optimization procedure is initialized with rejection sampling as described above for gradient approaches. However, to keep the optimized parameters within the feasible region, a subroutine is used to project infeasible offsprings onto the feasibility boundary using binary search (SI Sec.~\ref{sisubsubsec:genetic_acquisition}). In addition to guaranteeing that the optimized parameters satisfy the constraint function, this approach also ensures sampling of parameters close to the feasibility boundary (Fig.~\ref{fig:acq_demonstration}). We also use modified prior kernel precisions as described in the previous section.

Despite having implemented constraints in \gryffin independently, when preparing this manuscript we realized that known constraints in \dragonfly\cite{Kandasamy:2020} have been implemented following a very similar strategy. However, rather than projecting infeasible offspring solution on the feasibility boundary, \dragonfly relies solely on rejection sampling.

A practical advantage of the genetic optimization strategy is its favourable computational scaling compared to its gradient-based counterpart. We conducted benchmarks in which the time needed by \gryffin to optimize the acquisition function was measured for different numbers of past observations and parameter space dimensions (SI Sec.~\ref{sisubsubsec:acq_scaling}). In fact, acquisition optimization is the most costly task in most Bayesian optimization algorithms, including \gryffin. The genetic strategy provided a speedup of approximately $5\times$ over Adam when the number of observations was varied, and of approximately $2.5\times$ when the dimensionality of the optimization domain was varied. The better time complexity of the zeroth-order approach is primarily due to derivatives of $\alpha(\bm{x})$ not having to be computed. In fact, our Adam implementation currently computes derivatives numerically. Future work will focus on deriving the analytical gradients for \gryffin's acquisition function, or taking advantage of automatic differentiation schemes, such that this gap might reduce or disappear in future versions of the code.

\subsection{User interface} \label{subsec:method_interface}

With user-friendliness and flexibility in mind, we extended \gryffin's Python interface such that it can take a user-defined constraint function among its inputs. As an example, the following code snippet shows an instantiation of \gryffin where the sum of two volumes is constrained.

\begin{minted}[mathescape,autogobble,numbersep=5pt,
               frame=lines,framesep=2mm, fontsize=\small]{python}
from gryffin import Gryffin

# set path to optimization configuration
config_file = "config.json"

def my_constraints(params):
    vol0 = params["volume_0"]
    vol1 = params["volume_1"]
    if vol0 + vol1 > 50:
        return False
    else:
        return True
        
gryffin = Gryffin(config_file=config_file, 
                  known_constraints=my_constraints)
\end{minted}

This example could be that of a minimization over two parameters, $\bm{x} = (x_1, x_2) \in [0,1]^2$, subject to the constraint that the sum of $x_1$ and $x_2$ does not exceed some upper bound $b$,
     
\begin{align*} 
     \bm{x}^* & = \argmin{f(\bm{x})} \,. \\ 
     & \text{s.t.} \;\; x_1 + x_2 \leq b \nonumber 
\end{align*}

\subsection{Current limitations} \label{subsec:method_limitations}
We note two limitations of our approach to known constraints. First, in continuous spaces, the current implementation handles only inequality constraints. Equality constraints, such as $x_1 + x_2 = b$, effectively change the dimensionality of the optimization problem and can be tackled via a re-definition of the optimization domain or analytical transforms. For instance, for the two-dimensional example above, the constraint can be satisfied simply by optimizing over $x_1$, while setting $x_2 = b - x_1$. 

Second, because of the rejection sampling procedure used (Sec.~\ref{subsec:const_adam}), the cost of acquisition optimization is inversely proportional to the fraction of the optimization domain that is feasible. The smaller the feasible volume fraction is, the more samples need to be drawn uniformly at random before reaching the predefined number of samples to be collected for acquisition optimization. It thus follows that the feasible region, as determined by the user-defined constraint function, should not be exceedingly small (e.g., less than 1\% of the optimization domain), as in that case \gryffin's acquisition optimization would become inefficient. In practice we find that when only a tiny fraction of the overall optimization domain is considered feasible, it is because of a sub-optimal definition of the optimization domain, which can be solved by re-framing the problem. A related issue may arise if the constraints define disconnected feasible regions with vastly different volumes. In this scenario, given uniform sampling, little or no samples may be drawn from the smaller region. However, this too is an edge case that is unlikely to arise in practical applications, as the definition of completely separate optimization regions with vastly different scales tends to imply a scenario where multiple separate optimizations would be a more suitable solution.

Future work will focus on overcoming these two challenges. For instance, deep learning techniques like invertible neural networks\cite{Behrmann:2019} may be used to learn a mapping from the unconstrained to the constrained domain, such that the optimization algorithm would be free to operate on the hypercube while the proposed experimental conditions would satisfy the constraints.


\section{Results and Discussion} \label{sec:results}

\begin{figure*}[htb]
    \centering
    \includegraphics[width=1.0\columnwidth]{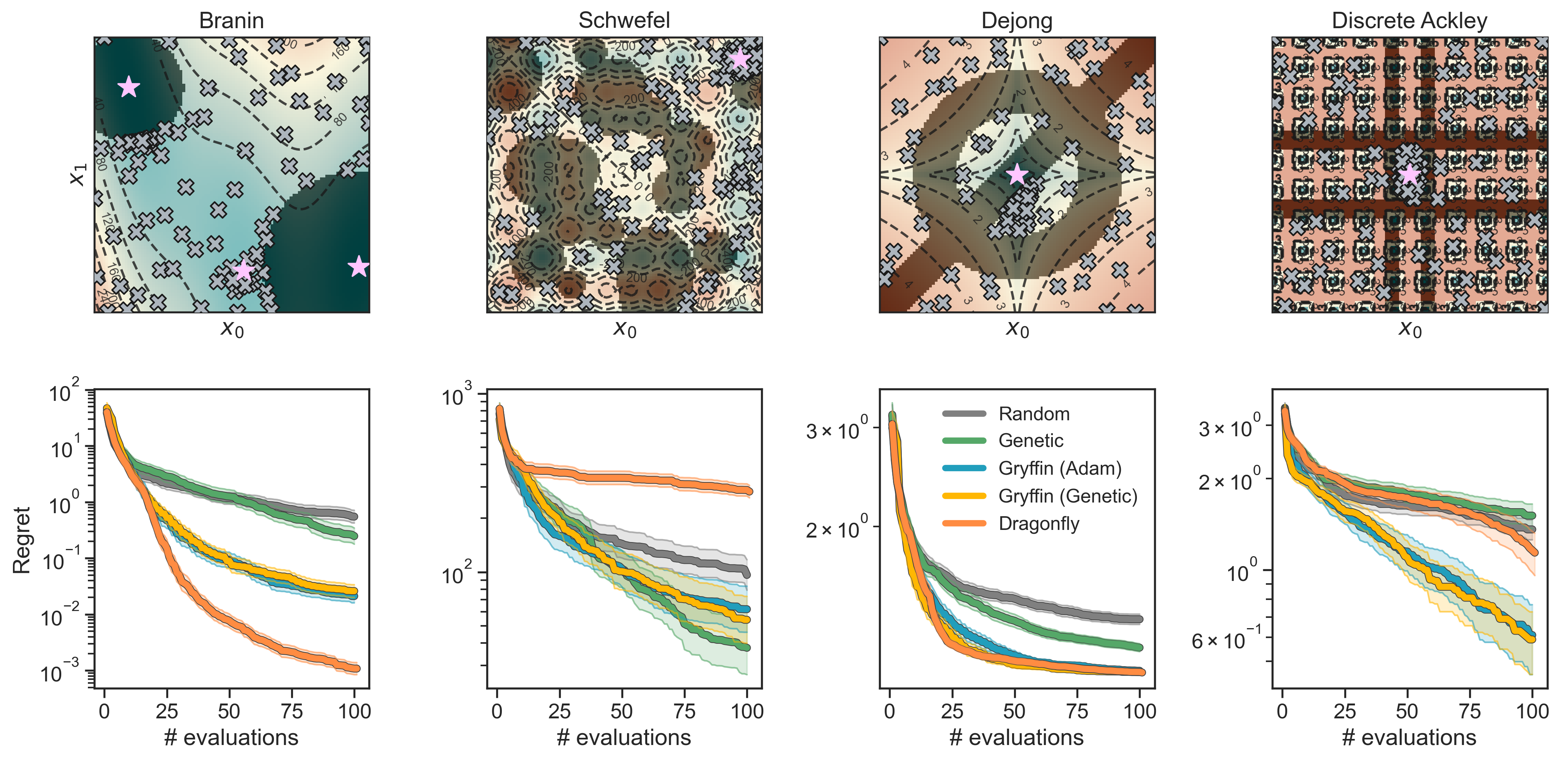}
    \caption{Constrained optimization benchmarks on analytical functions with continuous parameters. The upper row shows contour plots of the surfaces with constrained regions darkly shaded. Gray crosses show sample observation locations using the \textit{Gryffin (Genetic)} strategy and purple stars denote the location(s) of unconstrained global optima. The bottom row show optimization traces for each strategy. Shaded regions around the solid trace represent 95\% confidence intervals.} 
    \label{fig:synth_continuous}
\end{figure*}

In this section we test the ability of \gryffin to handle arbitrary known constraints. First, we show how \gryffin can efficiently optimize continuous and discrete (ordered) benchmark functions with a diverse set of constraints. Then, we demonstrate the practical utility of handling known constraints on two relevant chemistry examples: the optimization of the synthesis of \textit{o}-xylenyl Buckminsterfullerene adducts under constrained flow conditions, and the design of redox active molecules for flow batteries under synthetic accessibility constraints.

In addition to testing \gryffin---with both gradient and evolutionary based-acquisition optimization strategies, referred to as \textit{Gryffin (Adam/Hill)} and \textit{Gryffin (Genetic)}, respectively---we also test three other constrained optimization strategies amenable to experiment planning. Specifically, we use random search (\textit{Random}), a genetic algorithm (\textit{Genetic}), and Bayesian optimization with a Gaussian process surrogate model (\textit{Dragonfly}). \textit{Genetic} is the same algorithm developed for constrained acquisition function optimization, but it is here employed directly to optimize the objective function (Sec.~\ref{subsec:const_genetic}). \textit{Dragonfly} uses the \dragonfly package~\cite{Kandasamy:2020} for optimization. Similar to \gryffin, \dragonfly allows for the specification of arbitrarily complex constraints \textit{via} a Python function, which is not the case for most other Bayesian optimization packages. \textit{Dragonfly} was not employed in the two chemistry examples, due to an implementation incompatibility with using both constraints and the multi-objective optimization strategy required by the applications.


\subsection{Analytical benchmarks} \label{sec:benchmark}

\begin{figure*}[!ht]
    \centering
    \includegraphics[width=1.0\columnwidth]{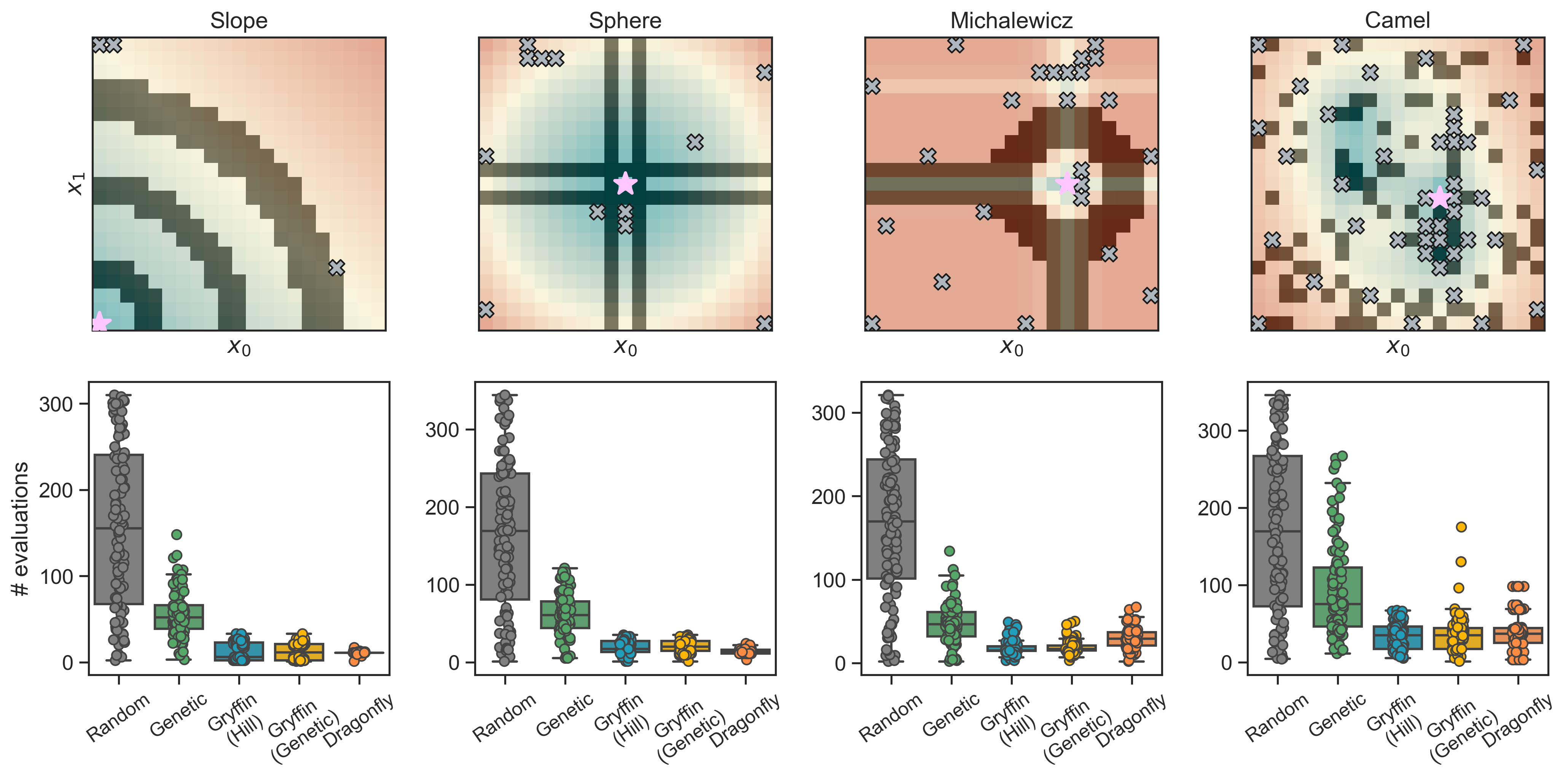}
    \caption{Constrained optimization benchmarks on analytical functions with discrete, ordered parameters. The upper row shows heatmaps of the discrete optimization domain. Shaded regions indicate infeasible regions, where constraints have been applied. The locations of the global optima are indicated by purple stars. Gray crosses show parameters locations that have been probed in a sample optimization run using \textit{Gryffin (Genetic)} before the optimum being identified. The bottom row shows, as superimposed box-and-whisker and swarm plots, the distributions of the number of evaluations needed to identify the global optimum for each optimization strategy (numerical values can be found in Table~\ref{sitab:cat_synth_opt_performances}).} 
    \label{fig:synth_categorical}
\end{figure*}

Here, we use eight analytical functions (four continuous and four discrete) to test the ability of \gryffin to perform sample-efficient optimization while satisfying a diverse set of user-defined constraints. More specifically, we consider the following four two-dimensional continuous surfaces, as implemented in the \textsc{Olympus} package~\cite{hase_olympus_2021}: \textit{Branin}, \textit{Schwefel}, \textit{Dejong}, and \textit{Discrete Ackley} (Fig.~\ref{fig:synth_continuous}). While \textit{Branin} has three degenerate global minima, we apply constraints such that only one minimum is present in the feasible region. In \textit{Schwefel}, a highly unstructured constraint function where the global minimum is close to the boundary of the infeasible region. In \textit{Dejong}, we use a set of constraints that result in a non-compact optimization domain where the global minimum cannot be reached, and there is an infinite number of feasible minima along the feasibility boundary. \textit{Discrete Ackley} is a discretized version of the Ackley function, which is an example of an extremely rugged surface. The four, two-dimensional discrete surfaces considered are: \textit{Slope}, \textit{Sphere}, \textit{Michalewicz}, and \textit{Camel} (Fig.~\ref{fig:synth_categorical}). Each variable in these discrete surfaces is comprised of integer numbers from $0$ to $20$, for a design space of $441$ options in total.

Results of constrained optimization experiments on continuous surfaces are shown in Fig.~\ref{fig:synth_continuous}. These results were obtained by performing 100 repeated optimizations for each of the five strategies considered, while allowing a maximum of 100 objective function evaluations. Optimization performance is assessed using the distance between the best function value found at every iteration of the optimization campaign and the global optimum, a metric known as \textit{regret}, $r$. The regret after $k$ optimization iterations is

\begin{align} \label{eq:regret}
 r_k = \vert f(\bm{x}^*)- f(\bm{x}_k^+)   \vert \,.
\end{align}

$\bm{x}_k^+$ are the parameters associated with best objective value observed in the optimization campaign after $k$ iterations, sometimes referred to as the \textit{incumbent} point, i.e. for a minimization problem $\bm{x}_k^+ = \argmin_{\bm{x} \in \mathcal{D}_k} f(\bm{x})$, where $\mathcal{D}_k$ is the current dataset of observations. $\bm{x}^*$ are the parameters associated with the global optimum of the function. Sample-efficient algorithms should find parameters that return better (in this case lower) function values with fewer objective evaluations.

All optimization strategies tested obeyed the known constraints (Fig.~\ref{fig:synth_continuous}, top row). We observe that the combination of analytical function and user-defined constraint have great influence on the relative optimization performance of the considered strategies. On all functions, the performance of \gryffin was insensitive to the acquisition optimization strategy, with \textit{Gryffin (Adam)} and \textit{Gryffin (Genetic)} performing equally well. On smooth functions, such as \textit{Branin} and \textit{Dejong}, \textit{Dragonfly} displayed strong performance, given that its Gaussian process surrogate model approximates these functions extremely well. In addition, \dragonfly appears to propose parameter points that are closer to previous observations than \gryffin does, showcasing a more exploitative character, which in this case is beneficial to achieve marginal gains in performance when in very close proximity to the optimum (see SI Sec. \ref{sisubsubsec:gryffin_dragonfly_comparison} for additional details). Although their performance is slightly worse than \textit{Dragonfly}, \textit{Gryffin} strategies still displayed improved performance over \textit{Random} and \textit{Genetic} on \textit{Branin}, and comparable perforamance to \textit{Dragonfly} on \textit{Dejong}. \textit{Gryffin} strategies showed improved performance compared to \textit{Dragonfly} on \textit{Schwefel} and \textit{DiscreteAckley}. These observations are consistent with previous work\cite{hase_phoenics_2018} where it was observed that, due to the different surrogate models used, \gryffin returned better performance on functions with discontinuous character, while Gaussian process-based approaches better performed on smooth surfaces. On \textit{Schwefel}, our \textit{Genetic} strategy also showed improved performance over that of \textit{Random} and \textit{Dragonfly}, comparable with that of \textit{Gryffin}. Finally, note that in Fig.~\ref{fig:synth_continuous} differences between strategies are exaggerated by the use of a logarithmic scale, used to highlight statistically significant differences. We also report the same results using a linear scale, which de-emphasizes significant yet marginal differences in regret values (e.g. between \textit{Dragonfly} and \textit{Gryffin} on the \textit{Branin} function) in SI Fig.~\ref{sifig:continuous_benchmarks_linear}.

Results of constrained optimization experiments on discrete synthetic surfaces are shown in Fig.~\ref{fig:synth_categorical}. These results are also based on $100$ repeated optimizations, each initialized with a different random seed. Optimizations were allowed to continue until the global optimum was found. As a measure of performance, the number of objective function evaluations required to identify the optimum was used. As such, a more efficient algorithm should identify the optimum with fewer function evaluations on average.

Here too, all strategies tested correctly obeyed each constraint function (Fig.~\ref{fig:synth_categorical}, top row), and Bayesian optimization algorithms (\textit{Gryffin} and \textit{Dragonfly}) outperformed \textit{Random} and \textit{Genetic} on all four benchmarks (Fig.~\ref{fig:synth_categorical}, bottom row). \textit{Random} needed to evaluate about half of the feasible space before identifying the optimum. \textit{Genetic} required significantly fewer evaluations, generally between a half and a third of those required by \textit{Random}. On discrete surfaces, the approaches based on \gryffin and \dragonfly returned equal performance overall, with no considerable differences (Table~\ref{sitab:cat_synth_opt_performances}). As it was observed for continuous surfaces, the performance of \textit{Gryffin (Adam)} and \textit{Gryffin (Genetic)} was equivalent.


\subsection{Process-constrained optimization of \textit{o}-xylenyl C\textsubscript{60} adducts synthesis} \label{subsec:fullerene_application}

\begin{figure*}[!p]
    \centering
    \includegraphics[width=1.0\textwidth]{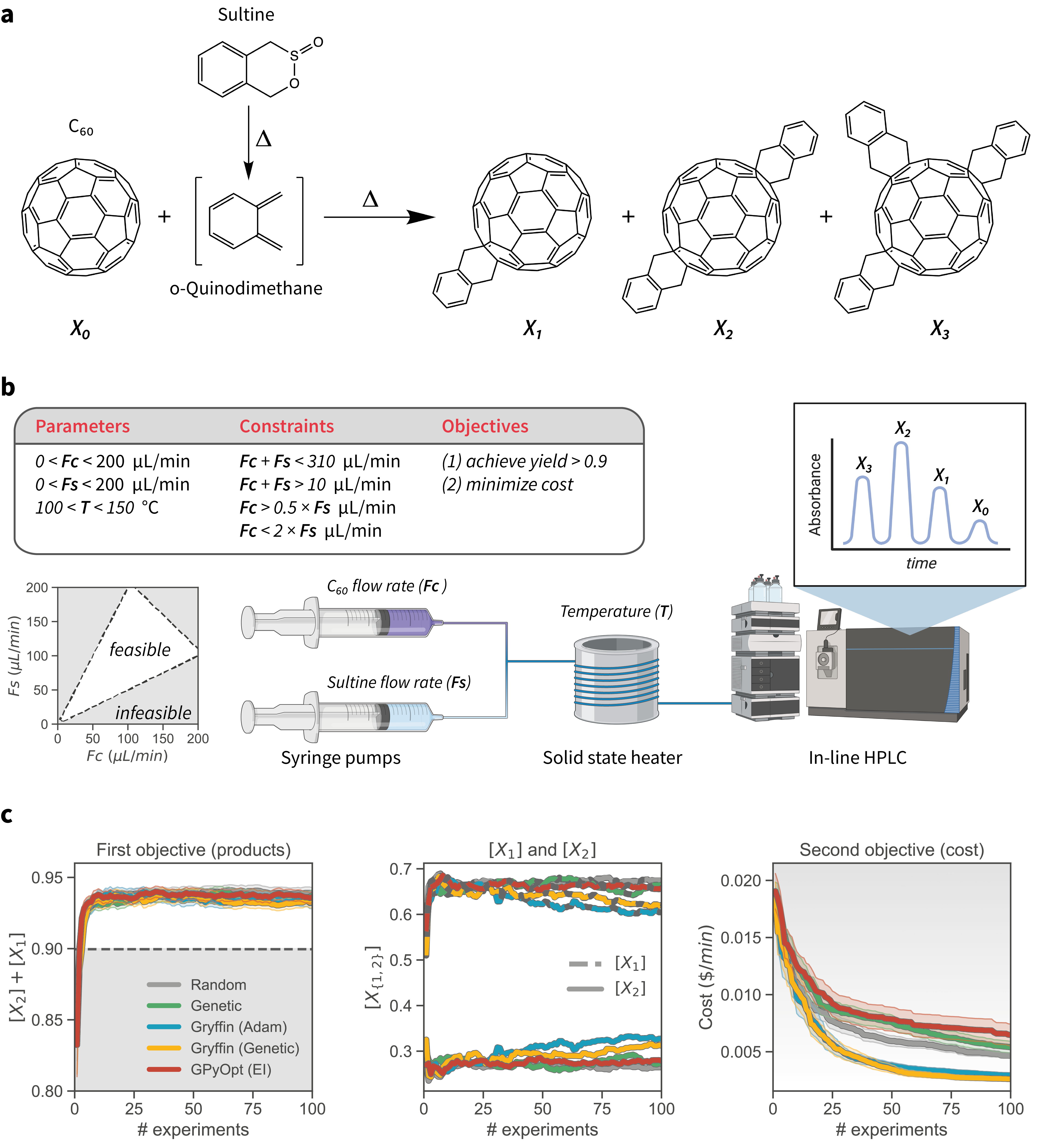}
    \caption{Experimental setup and results of the process-constrained synthesis of \textit{o}-xylenyl C\textsubscript{60} adducts. (a) Synthesis of \textit{o}-xylenyl C\textsubscript{60} adducts. 1,4-dihydro-2,3-benzoxathiin 3-oxide (sultine) is converted to \textit{o}-quinidimethane (\textit{o}QDM) \textit{in situ}, which then reacts with C\textsubscript{60} by Diels-Alder cycloaddition to form the \textit{o}-xylenyl adducts.  (b) Schematic of the single-phase capillary-based flow reactor (as reported by Walker \textit{et al.}\cite{walker_tuning_2017}), along with the optimization parameters, constraints, and objectives. The C\textsubscript{60} and sultine flow rates are modulated by syringe pumps, and the reaction temperature is controlled by a solid state heater. Online HPLC and absorption spectroscopy is used for analysis. (c) Results of the constrained optimization experiments. Plots show the mean and 95\% confidence interval of the best objective values found after varying numbers of experiments, for each optimization algorithm studied. The background color indicates the desirability of the objective values, with grey being less, and white more desirable. For the first objective (left-hand side), the grey region below the value of $0.9$ indicates values for the objective that do not satisfy the optimization goal that was set.}
    \label{fig:fullerenes_multipanel}
\end{figure*}

Compared to their inorganic counterparts, organic solar cells have the advantage of being flexible, lightweight and easily  fabricated.\cite{Chen:2019,Hong:2020} Bulk heterojunction polymer-fullerene cells are an example of such devices, in which the photoactive layer is composed of a blend of polymeric donor material and fullerene derivative acceptor molecules.\cite{Brabec:2010} The fullerene acceptor is oftentimes functionalized to tune its optoelectronic properties. In particular, mono and bis \textit{o}-xylenyl adducts of Buckminsterfullerene are acceptor molecules which have received much attention in this regard. On the other hand, higher-order adducts are avoided as they have a detrimental impact on the power conversion efficiencies of the resulting device.~\cite{kang_effect_2013} Therefore, synthetic protocols that provide accurate control over the degree of C\textsubscript{60} funcionalization are of primary interests for the manufacturing of these organic photovoltaic devices.

In this example application of constrained Bayesian optimization, we consider the reaction of C\textsubscript{60} with an excess of sultine, the cyclic ester of a hydroxy sulfinic acid, to form first, second, and third order \textit{o}-xylenyl adducts of Buckminsterfullerene (Fig.~\ref{fig:fullerenes_multipanel}a). This application is based on previous work by Walker \textit{et  al.},\cite{walker_tuning_2017} who reported the optimization of this reaction with an automated, single-phase capillary-based flow reactor. Because this reaction is cascadic, a mixture of products is always present, with the relative amounts of each species depending primarily on reaction kinetics. The overall goal of the optimization is thus to tune reaction conditions such that the yield of first- and second-order adducts is maximized and reaches at least 90\%, while the cost of reagents is minimized. We estimated reagents cost by considering the retail price listed by a chemical supplier (\ref{sisubsubsec:fullerene_cost_estimation}). In effect, we derive a simple estimate of per-minute operation cost of the flow reactor, which is to be minimized as the second optimization objective. 

The controllable reaction conditions are the temperature ($T$), the flow rate of C\textsubscript{60} ($F_{\text{C}}$), and the flow rate of sultine ($F_{\text{S}}$), which determine the chemical composition of the reaction mixture and are regulated by syringe pumps (Fig.~\ref{fig:fullerenes_multipanel}b). $T$ is allowed to be set anywhere between $100$ and $150$\textdegree C, and flow rates can be varied between $0$ and $200\ \mu\text{L}/\text{min}$. However, the values of the flow rates are constrained. First, the total flow rate cannot exceed $310 \; \mu \text{L}/\text{min}$, or be below $10 \; \mu \text{L}/\text{min}$. Second, $F_{\text{S}}$ cannot be more than twice $F_{\text{C}}$, and vice-versa $F_{\text{C}}$ cannot be more than twice $F_{\text{S}}$. More formally, the inequality constraints can be defined as $10 < F_{\text{C}}+F_\text{S} < 310 \; \mu \text{L}/\text{min}$, $ F_{\text{C}} < 2 F_\text{S}$ and $ F_\text{S} < 2 F_{\text{C}}$.

The relative concentrations of each adduct---$[X_1]$, $[X_2]$, and $[X_3]$ for the mono, bis, and tris adduct, respectively---are measured via online high performance liquid chromatography (HPLC) and absorption spectroscopy. Following from the above discussion, we would like to maximize the yield of $[X_1]$ and $[X_2]$, which are the adducts with desirable optoelectronic properties. However, we would also like to reduce the overall cost of the raw materials. These multiple goals are captured via the use of \chimera as a hierarchical scalarizing function for multi-objective optimization.\cite{hase_chimera_2018} Specifically, we set the first objective as the maximization of $[X_1] + [X_2]$, with a minimum target goal of $0.9$, and the minimization of reagents cost as the secondary goal. Effectively, this setup corresponds to the minimization of cost under the constraint that $[X_1] + [X_2] \geq 0.9$.

To perform numerous repeated optimizations with different algorithms, as well as collect associated statistics, we constructed a deep learning model to simulate the above-described experiment. In particular, we trained a Bayesian neural network based on the data provided by Walker \textit{et  al.},\cite{walker_tuning_2017} which learns to map reaction conditions to a distribution of experimental outcomes (SI Sec.~\ref{sisubsubsec:fullerene_emulator}). The measurement is thus stochastic, as expected experimentally. This emulator takes $T$, $F_{\text{C}}$, and $F_\text{S}$ as input, and returns mole fractions of the un- ($X_0$), singly- ($X_1$), doubly- ($X_2$), and triply-functionalized ($X_3$) C\textsubscript{60}. The model displayed excellent interpolation performance across the parameter space for each adduct type (Pearson coefficient of $0.93-0.96$ on the test sets). This enabled us to perform many repeated optimizations thanks to rapid and accurate simulated experiments.

Results of the constrained optimization experiments are shown in Fig.~\ref{fig:fullerenes_multipanel}c. Optimization traces show the objective values associated with the best scalarized merit value. The constrained \gryffin strategies obeyed all flow rate constraints defined. All strategies rapidly identified reaction conditions producing high yields. Upon satisfying the first objective, \chimera guided the optimization algorithms towards lowering the cost of the reaction. \gryffin achieved cheaper protocols faster than the other strategies tested.
In SI Sec.~\ref{sisubsubsec:fullerene_cost_estimation} we examine the best reactions conditions identified by each optimization strategy after 100 experiments. We find that for the majority of cases, strategies decreased the C\textsubscript{60} flow rate (the more expensive reagent) to minimize the cost of the experiment, while simultaneously decreasing the sultine flow rate to maintain a stoichiometric ratio close to one and preserve the high ($\geq 0.9$ combined mole fraction) yield of $X_0$ and $X_1$ adducts.
In principle, the yield is allowed to degrade with an accompanying decrease (improvement) in cost. However, we did not see the trade-off between yield and cost being required within the experimental budget of $100$ experiments. In this case study we used \gpyopt, as opposed to \dragonfly, as a representative of an established Bayesian optimization algorithm, because its implementation was compatible with the requirements of this specific scenario. That is, the presence of known constraints and the use of an external scalarizing function for multi-objective optimization.
Use of \chimera as an external scalarizing function requires updating the entire optimization history at each iteration. Thus, the optimization library must be used via an \textit{ask-tell} interface. Although \dragonfly does have such an interface, to the best of our knowledge it does not yet support constrained and multi-objective optimizations with external scalarizing functions.
In this instance, \gpyopt with expected improvement as the acquisition function performed similarly to \textit{Genetic} and \textit{Random} on average.

\subsection{Design of redox-active materials for flow batteries with synthetic accessibility constraints} \label{subsec:battery_application}

\begin{figure*}[!ht]
    \centering
    \includegraphics[width=1.0\columnwidth]{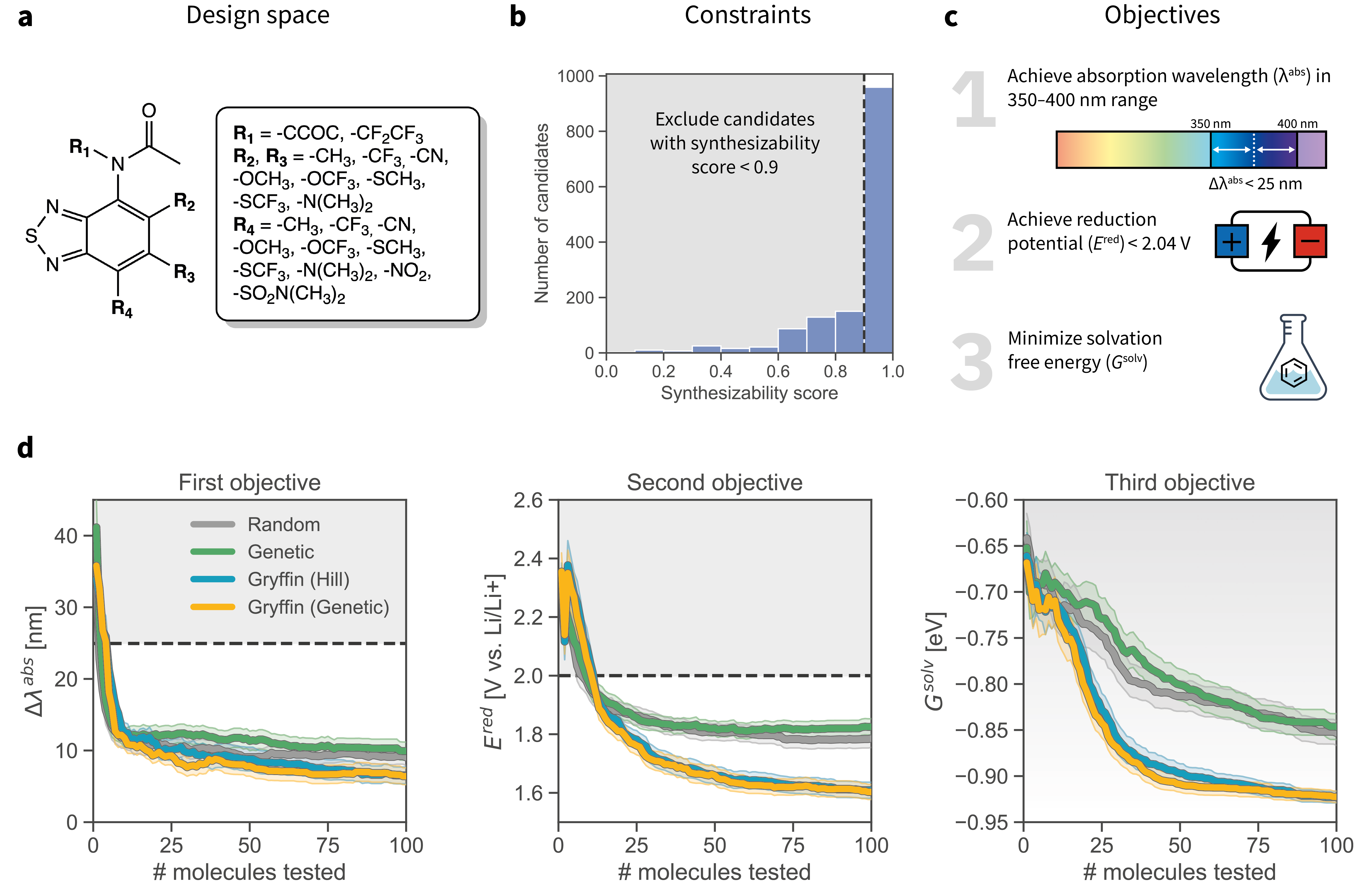}
    \caption{Optimization setup  and  results for the design of redox-active materials. (a) Markush structure of the benzothiadiazole scaffold and all subtituents considered. The entire design space consists of 1408 candidates (2 R\textsubscript{1} $\times$ 8 R\textsubscript{2} $\times$ 8 R\textsubscript{3} $\times$ 11 R\textsubscript{4} options). (b) Synthetic constraints applied to the optimization. The \textit{RAscore}\cite{thakkar_retrosynthetic_2021} was used to score synthesizability. The histogram shows the distribution of synthesizability scores for all $1408$ candidates. We constrain the optimization to those candidates with an \textit{RAscore} $> 0.9$. (c) Objectives of the molecular optimization. (d) Results of the constrained optimization experiments. Grey shaded regions indicate objective values failing to achieve the desired objectives. Traces depict the objective values corresponding to the best achieved merit at each iteration, where error bars represent 95\% confidence intervals.} 
    \label{fig:batteries_multipanel}
\end{figure*}

Long-duration, stationary energy storage devices are needed to handle the rapid growth in intermittent energy sources.~\cite{gur_review_2018} Toward this goal, redox flow batteries offer a potentially promising solution. \cite{lin_redox-flow_2016,leung_recent_2017,ye_redox_2017,lourenssen_vanadium_2019,kwabi_electrolyte_2020} Second-generation non-aqueous redox flow batteries are based on organic solvent as opposed to water, which enables a much larger electrochemical window and higher energy density, the potential to increase the working temperature window, as well as the use of cheap, earth abundant redox-active materials.~\cite{gong_nonaqueous_2015} Nevertheless, the multiobjective design of redox-active materials poses a significant challenge in materials science. Recently, the anolyte redox material 2,1,3-benzothiadiazole was shown to have low redox potential, high stability and promising electrochemial cycling performance.~\cite{duan_wine-dark_2017,huang_substituted_2018,zhang_comparing_2018} Furthermore, derivatization of benzothiadiazole enabled self-reporting of battery health by fluorescence emission.~\cite{robertson_fluorescence-enabled_2020}

In this application we demonstrate our constrained optimization approach for the multi-objective design of retro-synthetically accessible redoxmer molecules. We utilize a previously reported dataset that comprises $1408$ benzothiadiazole derivatives (Fig.~\ref{fig:batteries_multipanel}a) for which the reduction potential $E^{\text{red}}$, solvation free energy $G^{\text{solv}}$, and maximum absorption wavelength $\lambda^{\text{abs}}$ were computed with DFT.\cite{agarwal_discovery_2021} This example application thus simulates a DFT-based computational screen\cite{cruz_new_2020,bachman_investigation_2014,assary_reduction_2014} that aims to identify redoxmer candidates with self-reporting features and a high probability of synthetic accessibility.

To impose synthetic accessibility constraints we used \textit{RAscore}~\cite{thakkar_retrosynthetic_2021}, a recently-reported synthetic accessibility score based on the computer-aided synthesis planner \textit{AiZythFinder}~\cite{genheden_aizynthfinder_2020}. \textit{RAscore} predicts the probability of \textit{AiZythFinder} being able to identify a synthetic route for a target organic molecule. While other measures of synthetic accessibility are available,\cite{ertl_estimation_2009,coley_scscore_2018,vorsilak_syba_2020,seifrid_routescore_2021} \textit{RAscore} was chosen for its performance and intuitive interpretation.
In our experiments, we constrained the optimization to candidates with \textit{RAscore} $> 0.9$ (Fig.~\ref{fig:batteries_multipanel}b, additional details in SI Sec.~\ref{sisubsubsec:battery_synthetic_access}). This synthetic accessibility constraint reduces the design space of feasible candidates to a total of $959$ molecules.

We aimed at optimizing three objectives concurrently: the absorption wavelength $\lambda^{\text{abs}}$, the reduction potential $E^{\text{red}}$, and the solvation free energy $G^{\text{solv}}$ of the candidates (Fig.~\ref{fig:batteries_multipanel}c). Specifically, we aimed at identifying molecules that would (i) absorb in the 350$-$400 nm range, (ii) improve upon the reduction potential of the base scaffold ($2.04$ V against a Li/Li$^+$ reference electrode; SI Sec.~\ref{sisubsubsec:dft_reduction_constraint}), and (iii) provide the lowest possible solvation free energy, here used as a proxy for solubility. The hierarchical scalarizing function \chimera\cite{hase_chimera_2018} was used to guide the optimization algorithms toward achieving these desired properties.

Results of the constrained optimization experiments are displayed in Fig.~\ref{fig:batteries_multipanel}d. Each optimization strategy was given a budget of $100$ molecules to be tested; the experiments were repeated $100$ times and initialized with different random seeds. Optimization traces show the values of $\Delta \lambda^{\text{abs}}$, $E^{\text{red}}$, and $G^{\text{solv}}$ associated  with  the  best  molecule identified among those tested. In these tests, the dynamic version of \gryffin was used.\cite{hase_gryffin_2021} This approach can take advantage of physicochemical descriptors in the search for optimal molecules. In this case, \gryffin was provided with several descriptors for each R group in the molecule: number of heavy atoms, number of hetero atoms, molecular weight, geometric diameter, polar surface area, polarizability, fraction of sp$^3$ carbons (SI Sec.~\ref{sisubsubsec:battery_descriptors}).~\cite{moriwaki_mordred_2018}

Constrained \gryffin avoided redoxmer candidates predicted to be retro-synthetically inaccessible. All optimization strategies rapidly identified candidates with $\lambda^{\text{abs}}$ between 350 and 400 nm, and whose $E^{\text{red}}$ is lower than that of the starting scaffold. However, \gryffin strategies managed to identify candidates with lower $G^{\text{solv}}$ faster than \textit{Random} and \textit{Genetic}. Consistent with previous tests, we did not observe a statistically different performance between \textit{Gryffin (Hill)} and \textit{Gryffin (Genetic)}. Overall, \textit{Gryffin (Hill)} and \textit{Gryffin (Genetic)} identified molecules with better $\lambda^{\text{abs}}$, $E^{\text{red}}$, and $G^{\text{solv}}$ properties than \textit{Random} or \textit{Genetic} after probing $100$ molecules. Furthermore, within this budget, \gryffin identified molecules with better properties more efficiently than the competing strategies tested. Without the aid of physicochemical descriptors, \gryffin's performance deteriorated, as expected,\cite{hase_gryffin_2021} yet it was still superior than that of \textit{Random} and \textit{Genetic} (SI Fig.~\ref{sifig:batteries_traces}).


\section{Conclusion} \label{sec:conclusion}

    In this work, we extended the capabilities of \gryffin to handle \textit{a priori} known, hard constraints on the parameter domain, a pragmatic requirement for the development of autonomous research platforms in chemistry. Known constraints constitute an important class of restrictions placed on optimization domains, and may reflect physical constraints or limitations in laboratory equipment capabilities. In addition, known constraints may provide a straightforward avenue to to inject prior knowledge or intuition into a chemistry optimization task.
    In \gryffin, we allow the user to define arbitrary known constraints via a flexible and intuitive Python interface. The constraints are then satisfied by constraining the optimization of its acquisition function.
    In all our benchmarks, \gryffin obeyed the (sometimes complex) constraints defined, showed superior performance to more traditional search strategies based on random search and evolutionary algorithms, and was competitive to state-of-the-art Bayesian optimization approaches.
    Finally, we demonstrated the practical utility of handling known constraints with \gryffin in two research scenarios relevant to chemistry. In the first, we demonstrated how to perform an efficient optimization for the synthesis of \textit{o}-xylenyl adducts of Buckminsterfullerene while subjecting the experimental protocol to process constraints. In the second, we showed how known constraints may be used to incorporate synthetic accessibility considerations in the design of redox active materials for non-aqueous flow batteries.
    It is our hope that simple, flexible, and scalable software tools for model-based optimization over constrained parameter spaces will enhance the applicability of data-driven optimization in chemistry and material science, and will contribute to the operationalization of self-driving laboratories.

	\section*{Data availability}
	
	An open-source implementation of \gryffin with known constraints is available on GitHub at \hyperlink{https://github.com/aspuru-guzik-group/gryffin}{https://github.com/aspuru-guzik-group/gryffin}, under an Apache 2.0 license. The data and scripts used to run the experiments and produce the plots in this paper are also available on GitHub at \hyperlink{https://github.com/aspuru-guzik-group/gryffin-known-constraints}{https://github.com/aspuru-guzik-group/gryffin-known-constraints}.
	

	\section*{Acknowledgments}

    The authors thank Dr. Martin Seifrid, Marta Skreta, Ben Macleod, Fraser Parlane, and Michael Elliott for valuable discussions.
    R.J.H. gratefully acknowledges the Natural Sciences and Engineering Research Council of Canada (NSERC) for provision of the Postgraduate Scholarships-Doctoral Program (PGSD3-534584-2019). 
    M.A. was supported by a Postdoctoral Fellowship of the Vector Institute.
    A.A.-G. acknowledges support from the Canada 150 Research Chairs program and CIFAR, as well as the generous support of Anders G. Fr\"oseth.  
    This work relates to the Department of Navy award (N00014-18-S-B-001) issued by the office of Naval Research. The United States Government has a royalty-free license throughout the world in all copyrightable material contained herein. Any opinions, findings, and conclusions or recommendations expressed in this material are those of the authors and do not necessarily reflect the views of the Office of Naval Research.
    Computations reported in this work were performed on the computing clusters of the Vector Institute and on the Niagara supercomputer at the SciNet HPC Consortium.~\cite{niagara1,niagara2} Resources used in preparing this research were provided, in part, by the Province of Ontario, the Government of Canada through CIFAR, and companies sponsoring the Vector Institute. SciNet is funded by the Canada Foundation for Innovation, the Government of Ontario, Ontario Research Fund - Research Excellence, and by the University of Toronto.

	\section*{Conflicts of interest}
	
	A.A.-G is the Chief Visionary Officer and founding member of Kebotix, Inc.


    \phantomsection\addcontentsline{toc}{section}{\refname}\putbib[main]
\end{bibunit}

\clearpage
\newpage

\begin{bibunit}[ieeetr]
	
	\setcounter{page}{1}
	\setcounter{section}{0}
	\setcounter{subsection}{0}
	\setcounter{figure}{0} 
	
        \renewcommand{\thesection}{S.\arabic{section}}
	\renewcommand{\thefigure}{S\arabic{figure}}
	\renewcommand{\thetable}{S\arabic{table}}
	
	\onecolumngrid
	\subsection*{\normalsize{Supplementary Information}{\\}{\vspace{6pt}}
	 			\large{Bayesian optimization with known experimental and design constraints for chemistry applications}{\\}{\vspace{6pt}}
			        \normalsize{\normalfont{Riley J. Hickman,$^{1,2,*}$ Matteo Aldeghi,$^{1,2,3,4,*}$ Florian H\"{a}se,$^{1,2,3,5}$ Al\'{a}n Aspuru-Guzik$^{1,2,3,6,7,8}$}}
			        {\\}{\vspace{6pt}}
			        \small{\normalfont{$^1$\textit{Chemical Physics Theory Group, Department of Chemistry, University of Toronto, Toronto, ON M5S 3H6, Canada}}}{\\}
			        \small{\normalfont{$^2$\textit{Department of Computer Science, University of Toronto, Toronto, ON M5S 3G4, Canada}}}{\\}
			        \small{\normalfont{$^3$\textit{Vector Institute for Artificial Intelligence, Toronto, ON M5S 1M1, Canada}}}{\\}
			        \small{\normalfont{$^4$\textit{Department of Chemical Engineering, Massachusetts Institute of Technology, Cambridge, MA 02139, United States}}}{\\}
			        \small{\normalfont{$^5$\textit{Department of Chemistry and Chemical Biology, Harvard University, Cambridge, MA 02138, United States}}}{\\}
			        \small{\normalfont{$^6$\textit{Department of Chemical Engineering \& Applied Chemistry, University of Toronto, Toronto, ON M5S 3E5, Canada}}}{\\}
			        \small{\normalfont{$^7$\textit{Department of Materials Science \& Engineering, University of Toronto, Toronto, ON M5S 3E4, Canada}}}{\\}
			        \small{\normalfont{$^8$\textit{Lebovic Fellow, Canadian Institute for Advanced Research, Toronto, ON M5G 1Z8, Canada}}}{\\}
			        \small{\normalfont{$^*$\textit{These authors contributed equally}}}
			    }
	{\vspace{18pt}}



\section{Constrained optimization of the acquisition function} \label{sisubsec:constrained_acquisition}

Minima of the acquisition function, $\alpha(\bm{x})$  (main text Eq.~\ref{eq:gryffin_acq}), determine the parameter space point $\bm{x}$ to be proposed for objective function measurement, thus, its optimization is an important subroutine in Bayesian optimization.
Acquisition function optimization in \gryffin consists of two main steps. First, a set of $N$ initial parameter space points are drawn from the domain $\mathcal{X}$ using rejection sampling according to constraint function $c(\bm{x})$ producing the set of initial samples $\mathcal{P}_{\text{init}}= \{ \bm{x}_i\}_{i=1}^{N}$, where $\bm{x}_i \sim \mathcal{X}; \; \text{s.t.}\ c(\bm{x}_i) \mapsto \text{feasible}$. These samples are then refined by one of several optimization strategies, returning a set of refined proposals $\mathcal{P}_{\text{ref}}= \{ \bm{x}_i\}_{i=1}^{N}$. In the following subsections, we detail the each acquisition function optimization strategy and compare their computational scaling. 

\subsection{Adam/Hill optimizer} \label{sisubsubsec:adam_acquisition}

For continuous parameters, the Adam optimizer~\cite{kingma_adam_2017} is used to optimize the acquisition function, with built in checks for whether the optimized samples obey the defined known constraints. For discrete and categorical parameter types, we use a $``$hill-climbing" strategy, in which each initial sample is randomly updated for a predefined number of iterations, with the candidate that has the best merit maintained and eventually returned. Algorithm~\ref{alg:acq_opt_hill} shows this strategy's basic pseudocode.

\begin{algorithm}
\caption{Constrained acquisition function optimization with Adam/Hill climbing}\label{alg:acq_opt_hill}
\KwData{initial samples $\mathcal{P}_{\text{init}}$, acquisition function $\alpha(\cdot)$, constraint function $c(\cdot)$, max iterations $i_{\text{max}}$}
\KwResult{refined samples $\mathcal{P}_{\text{ref}}$}
$\mathcal{P}_{\text{ref}} \gets \emptyset$ \;
\For{$\bm{x}_n$ in $\mathcal{P}_{\text{init}}$}{
    \For{$i$ in $i_{\text{max}}$}{
        $\bm{x}_n \gets \texttt{AdamStep}\left(\bm{x_n^{\text{cont}}}\right)$
        \Comment*[r]{continuous optimization step} 
        $\bm{x}_n \gets \texttt{HillStep}\left( \bm{x}_n^{\text{cat}} \right)$
        \Comment*[r]{categorical optimization step} 
        $\bm{x}_n \gets \texttt{HillStep}\left( \bm{x}_n^{\text{disc}} \right)$
        \Comment*[r]{discrete optimization step}
        \uIf{$c(\bm{x}_n) \mapsto \text{infeasible}$}{
            $\mathcal{P}_{\text{ref}} \overset{+}{\leftarrow} 
            \bm{x}_n^{i-1}$ \Comment*[r]{add feasible sample from previous iteration to refined samples} 
        }\uElse(){ 
            $\mathcal{P}_{\text{ref}} \overset{+}{\leftarrow} \bm{x}_n$ \Comment*[r]{add current sample to refined samples}
     }  
    }
}
\SetKwFunction{FMain}{HillStep}
\SetKwProg{Pn}{Function}{:}{\KwRet{$\bm{x}_n$}}
\Pn{\FMain{$\bm{x}_n$}}{
    $y_{\text{best}} \gets \alpha\left(\bm{x}_n\right)$\;
    \For{$\bm{z}_d$ in $\bm{x}_n^{\text{cat}}$}{
        $\bm{z}_d^{\text{new}} \gets \bm{z}_d^{\text{new}} \sim \mathbb{S}_d$ \Comment*[r]{sample from set of categorical options $\mathbb{S}_d$}
        $\bm{x}_n^{\text{new}}$ updated with $\bm{z}_d^{\text{new}}$ \;
        $y_{\text{new}} \gets \alpha\left( \bm{x}_n^{\text{new}}\right)$\;
        \uIf{$y_{\text{new}} < y_{\text{best}}$}{
            $y_{\text{best}} \gets y_{\text{new}}$ \;
            $\bm{x}_n \gets \bm{x}_n^{\text{new}}$ \;
        }
    }
}
\end{algorithm}

\subsection{Genetic optimizer} \label{sisubsubsec:genetic_acquisition}

Our genetic algorithm implementation is based on the \deap library~\cite{Fortin:2012,DeRainville:2012} and consists of crossover ($\mathcal{C}$), mutation ($\mathcal{M}$) and tournament selection ($\mathcal{S}$) operations. The size of the initial population is determined by the number of initial random samples, $N$. At each generation, offspring are chosen using tournament-style selection (tournament size of 3); elitism is applied to 5\% of the population. Each parent in the resulting population is given a chance to mate (using either uniform or two-point crossover, depending on the dimensionality of the parameter space, with a crossover probability of $0.5$), as well as a chance to mutate. We employ a custom mutation operator, $\mathcal{M}$ which can handle continuous, discrete and categorical parameter types. The mutation probability is set to $0.4$, with the probability of mutating each individual parameter set to $0.2$. The maximum number of generations is set to $10$, but the search is terminated early if the diversity of the population reaches a specified threshold. Specifically, if the population is concentrated in a small subvolume of parameter space, where each attribute does not span more than 10\% of the allowed range, then the search is terminated.
$\mathcal{M}$ consists of different mutations for each parameter type. For continuous parameters, a perturbation $x'$ is sampled from a Gaussian distribution with scale $0.1$, i.e. $\mathcal{M}(x) = x + x'$, $x' \sim \mathcal{N}(0, 0.1)$. The same is true for discrete parameter types, except $x'$ is first rounded to the nearest integer value. For categorical parameters, the mutated offspring are sampled randomly from the set of options for that parameter.

In order to constrain the genetic optimization procedure according to the known constraint function, a subroutine is used to project infeasible population offspring onto the feasibility boundary using a binary search procedure. After each application of $\mathcal{C}$ or $\mathcal{M}$ to parent $\bm{x}_p$, the feasibility of each resulting offspring $\bm{x}_o$ is tested with $c(\bm{x}_o)$, where $c$ is the user-defined constraints function. If $c(\bm{x}_o)$ returns \texttt{False},  i.e. the constraint is not satisfied and $\bm{x}_o$ is in the infeasible region, the following procedure is employed project $\bm{x}_o$ onto the boundary of the feasible region. For continuous parameters, we consider the midpoint $\bm{x}_m$ of the line segment bounded by $\bm{x}_p$ and $\bm{x}_o$. If $\bm{x}_m$ is feasible, then the parent is set to the midpoint; if it is infeasible, the offspring is set to the midpoint. This process is repeated until the distance between $\bm{x}_o$ and $\bm{x}_p$ is below a tolerance threshold. As a default, we use the criterion $\Vert \bm{x}_p - \bm{x}_o \Vert_{\infty} < 0.01$ to terminate the search, i.e., when we are guaranteed to be within a 1\% of relative distance from the feasibility boundary in all parameter dimensions. Throughout this procedure, $\bm{x}_p$ is guaranteed to always be feasible while $\bm{x}_o$ is always infeasible. Hence, once the search is terminated, the $\bm{x}_o$ is set equal to $\bm{x}_p$ and is returned. For discrete parameters, the same process is performed but the search is terminated when the closest point to the feasibility boundary is identified. Categorical parameters of infeasible offspring are instead simply reset to those of the feasible parent. When mixed continuous, discrete,   and categorical parameters are present, we (i) set the categorical parameters of $\bm{x}_o$ to those of its parent, to obtain $\bm{x}_o'$. If $\bm{x}_o'$ is feasible,   we return it, otherwise we (ii) perform the binary search procedure described above for the continuous and/or discrete parameters and obtain $\bm{x}_o''$. Then, we (iii) reset the categorical parameters of $\bm{x}_o''$ to their original values in $\bm{x}_o$,   obtaining $\bm{x}_o'''$. If $\bm{x}_o'''$ is feasible, we return it, otherwise we return $\bm{x}_o''$. Given this approach relies on binary searches, it has a favorable logarithmic scaling and adds little overhead to the optimization of the acquisition function. 

Algorithm~\ref{alg:acq_opt_genetic} shows the basic pseudocode for our implementation. We show pseudocode for our custom mutation function $\mathcal{M}$ (referred to as \texttt{Mutation}), but we omit definition of our subroutine which projects infeasible points to the feasible boundary for brevity. The function is referred to in Algorithm~\ref{alg:acq_opt_genetic} by \texttt{ProjectToFeasible}. We direct the interested reader to the source code of \gryffin for more details (\hyperlink{https://github.com/aspuru-guzik-group/gryffin}{https://github.com/aspuru-guzik-group/gryffin}).

\begin{algorithm}
\caption{Constrained acquisition function optimization with genetic algorithm}\label{alg:acq_opt_genetic}
\KwData{$\mathcal{P}_{\text{init}}$, $\alpha(\cdot)$, $c(\cdot)$, $i_{\text{max}}$, crossover operator $\mathcal{C}$ with prob $p_{\mathcal{C}}$, custom mutation operator $\mathcal{M}$ with prob $p_{\mathcal{M}}$ and independent prob $p_{\mathcal{M}}^{indep}$, tournament selection operator $\mathcal{S}$}
\KwResult{refined population $\mathcal{P}_{\text{ref}}$}
$\mathcal{P} \gets \mathcal{P}_{\text{init}}$; $f \gets \alpha \left(  \mathcal{P} \right)$ \;
\For{$i$ in $i_{\text{max}}$}{
    $\mathcal{O} \gets \mathcal{S} \left( P \right)$ \Comment*[r]{tournament selection of offspring $\mathcal{O}$ from population $\mathcal{P}$}
    \For{$\bm{x}^i_{parent,1}$, $\bm{x}^i_{parent,2}$ in $mating$ $pairs$}{
            \uIf{crossover sample $\sim \mathcal{U}\left(0, 1 \right) < p_{\mathcal{C}}$}{
                $\bm{x}^i_{child,1}, \bm{x}^i_{child,2} \gets \mathcal{C} \left( \bm{x}^i_{parent,1}, \bm{x}^i_{parent, 2} \right)$\;
                $\bm{x}^i_{child,1} \gets \texttt{ProjectToFeasible} \left( \bm{x}^i_{child,1}, \bm{x}^i_{parent,1}\right)$\;
                $\bm{x}^i_{child,2} \gets \texttt{ProjectToFeasible} \left( \bm{x}^i_{child,2}, \bm{x}^i_{parent,2}\right)$\;
            }
    }
    \For{$\bm{x}^i_{parent}$ in $\mathcal{O}$}{
        \uIf{mutant sample $\sim \mathcal{U}\left(0,1\right) < p_{\mathcal{M}}$}{
            $\bm{x}^i_{mutant} \gets \mathcal{M} \left( \bm{x}^i_{parent} \right)$\;
            $\bm{x}^i_{mutant} \gets \texttt{ProjectToFeasible} \left( \bm{x}^i_{mutant}, \bm{x}^i_{parent}\right)$\;
        }
    }
    $f \gets \alpha \left( \mathcal{O} \right)$ \Comment*[r]{evaluate the fitness $f$ of offspring $\mathcal{O}$} 
    $\mathcal{O} \overset{+}{\leftarrow} \mathcal{E}$ \Comment*[r]{add elites $\mathcal{E}$ to the offspring $\mathcal{O}$} 
    $\mathcal{P} \gets \mathcal{O}$ \Comment*[r]{set population $\mathcal{P}$ as the offspring $\mathcal{O}$ for next generation}
}
$\mathcal{P}_{\text{ref}} \gets \mathcal{P}$\;
\SetKwFunction{FMain}{Mutation}
\SetKwProg{Pn}{Function}{:}{\KwRet{$\bm{x}$}}
\Pn{\FMain{$\bm{x}$, $p_{\mathcal{M}}^{indep}$}}{
    \For{$x_d$ in $\bm{x}$}{
        \uIf{idependent mutation sample $\sim  \mathcal{U} \left(0,1\right) < p_{\mathcal{M}}^{indep}$}{
            \uIf{$x_d$ is continuous}{
                $x_d \gets x_d + x'$ \Comment*[r]{sample perturbation from unit Gaussian, i.e. $x' \sim \mathcal{N}(0,1)$}
              }
              \uElseIf{$x_d$ is discrete}{
                $x_d \gets x_d + \text{round}(x', \text{integer})$ \Comment*[r]{sample perturbation from unit Gaussian, i.e. $x' \sim \mathcal{N}(0,1)$}
              }
              \uElseIf{$x_d$ is categorical}{
                $x_d \gets x'$ \Comment*[r]{sample $x'$ from set of categorical options, i.e. $x' \sim \mathbb{S}_d$}
              }
        }
    }
}
\end{algorithm}

\subsection{Empirical time complexity of the Adam and genetic acquisition optimizers} \label{sisubsubsec:acq_scaling}

\begin{figure*}[!ht]
    \centering
    \includegraphics[width=1.0\columnwidth]{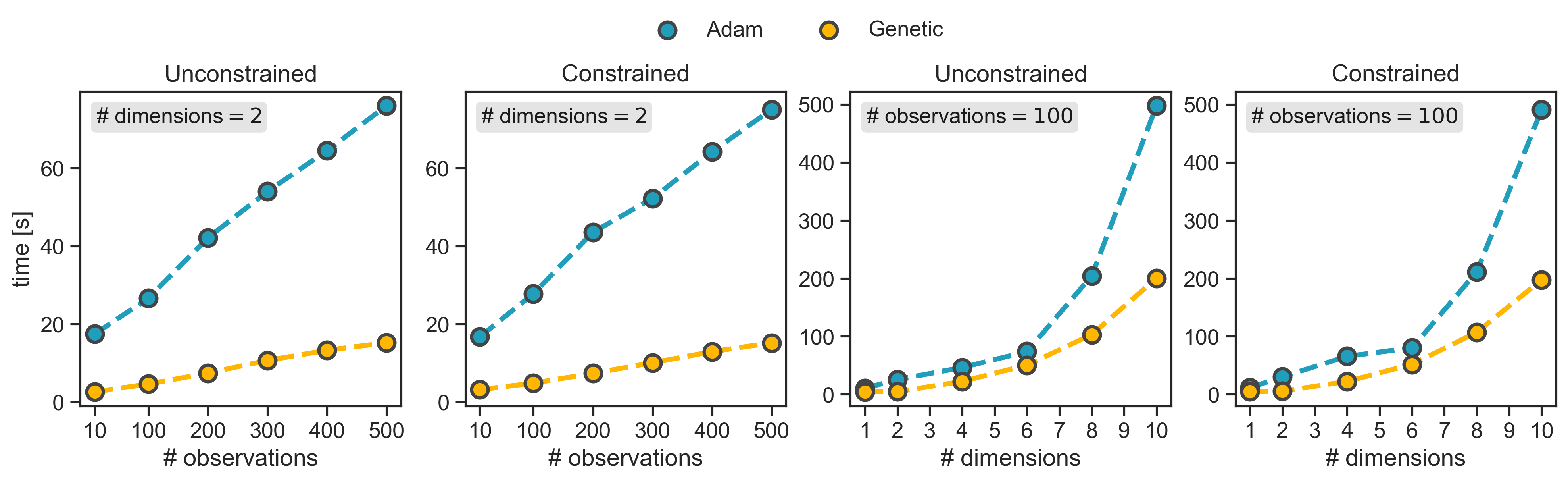}
    \caption{Empirical measurements of the time required by \gryffin to optimize its acquisition function. The Adam and Genetic optimization strategies were compared at varying number of past observations and optimization domain dimensions, with and without the presence of constraints.} 
    \label{sifig:acq_opt_scaling}
\end{figure*}

Computational  scaling experiments were carried out to compare the relative cost of the Adam and Genetic acquisition optimization strategies. The time taken to optimize the acquisition function was measured with increasing number of past observations while keeping the problem dimensionality constant, and with increasing number of parameter dimensions while keeping the number of observations constant. All parameters were continuous and in $[0,1]$. Tests were performed with and without optimization constraints. When relevant, the constraint used was $\sum_{i=1}^{d} x_i \leq 0.5d$  ,   where $d$ is dimension of the parameter space and $x_i$ are the individual elements of the $d$-dimensional parameter vector. That is, we assumed half of the optimization domain to be infeasible. Results of these tests are shown in Fig.~\ref{sifig:acq_opt_scaling}. Each datapoint was obtained as the average elapsed time for $60$ repeated acquisition function optimizations ($20$ for each $\lambda= \{-1,0,1\}$). No appreciable overhead was observed when constraints were present. When keeping the dimensionality constant, the Genetic strategy showed favourable scaling compared to Adam, being roughly 20\% as expensive as Adam after 100 observations. Similar results were observed in the experiments with a constant number of observations, where the optimization cost with the Genetic strategy took, on average, $\sim 60$\% less time than Adam.

\section{Constrained optimization of analytical functions} \label{sisubsec:synthetic_benchmarks}
In this section, we provide more details about the constrained analytical functions used for testing \gryffin's implementation and performance. 

\subsection{Benchmark functions and constraints used}

Our synthetic benchmark experiments consisted of four continuous and four discrete surfaces in two dimensions. The original implementations of the surfaces can be accessed via the \olympus package.~\cite{hase_olympus_2021} We used Python wrappers for each of the surfaces to implement constraints on the parameter space. While the full implementation is available on GitHub, code snippets are provided here as well to show a user may implement different constraint functions to be used by \gryffin. These constraint functions, called \texttt{is\_feasible()}, expect a dictionary,   \texttt{params}, containing the parameter values and evaluate their feasibility; \texttt{True} is returned for feasible, \texttt{False} for infeasible. Note that,   while here we report the definition of the analytical functions and the input domain typically used, \olympus normalizes the input domain to the unit hypercube for ease of use (i.e. all analytical function can be expected to be supported in $[0,1]^d$). As such, the constraint functions below assume each parameter to be normalized between zero and one.
\begin{itemize}

\item \textit{Branin} This surface is evaluated on the domain $x_1 \in [-5,10]$,   $x_2 \in [0,15]$, and has the form $f(\bm{x})=a(x_2 - bx_1^2 + cx_1 -r)^2 + s(1-t)\cos(x_1)+s$, with $a=1$, $b=5.1/4\pi^2$, $c=5/\pi$, $t=1/8\pi$. There are three degenerate global minima at $(x_1, x_2) = (-\pi,12.275)$, $(\pi,2.275)$ and $(9.42478, 2.475)$. Two of these minima were removed by the constraints defined below.

\begin{minted}[mathescape  ,  autogobble  ,  numbersep=5pt  ,   frame=lines  ,  framesep=2mm  ,   fontsize=\small]{python}
def is_feasible(params):
    x0 = params['x0']
    x1 = params['x1']
    y0 = (x0-0.12389382)**2 + (x1-0.81833333)**2
    y1 = (x0-0.961652)**2 + (x1-0.165)**2
    if y0 < 0.2**2 or y1 < 0.35**2:
        return False
    else:
        return True
\end{minted}

 \item \textit{Schwefel}: This surface is a complex optimization problem with many local minima. In $d$ dimensions, it is evaluated on the hypercube $x_i \in [-500, 500] \; \forall \; i=1,\ldots,d$ and is described by the expression $f(\bm{x}) = 418.9829 d - \sum_{i=1}^d x_i \sin \left( \sqrt{|x_i|}\right)$. The surface has a global optima at $\bm{x} = \left( 420.9687, \ldots, 420.9687 \right)$.

\begin{minted}[mathescape  ,  autogobble  ,  numbersep=5pt  ,   frame=lines  ,  framesep=2mm  ,   fontsize=\small]{python}
def is_feasible(params):
    np.random.seed(42)
    N = 20
    centers = [np.random.uniform(low=0.0, high=1.0, size=2) for i in range(N)]
    radii = [np.random.uniform(low=0.05, high=0.15, size=1) for i in range(N)]
    x0 = params['x0']
    x1 = params['x1']
    Xi = np.array([x0, x1])
    for c, r in zip(centers, radii):
        if np.linalg.norm(c - Xi) < r:
            return False
    return True
\end{minted}

 \item \textit{Dejong}: This surface generalizes a parabola to higher dimensions. It is convex and unimodal an evaluated on the $d$-dimensional hypercube $x_i \in [-5, 5],   \; \forall \; i=1, \ldots, d$. In two dimensions, this surface has a global minimum at $(x_0, x_1) = (0, 0)$ with $y=0$.
 
 \begin{minted}[mathescape  ,  autogobble  ,  numbersep=5pt  ,   frame=lines  ,  framesep=2mm  ,   fontsize=\small]{python}
def is_feasible(params):
    x0 = params['x0']
    x1 = params['x1']
    y = (x0-0.5)**2 + (x1-0.5)**2
    if np.abs(x0-x1) < 0.1:
        return False
    if 0.05 < y < 0.15:
        return False
    else:
        return True
\end{minted}

 \item \textit{DiscreteAckley}: This surface is the discrete analogue to the \textit{Ackley} function.

 \begin{minted}[mathescape  ,  autogobble  ,  numbersep=5pt  ,   frame=lines  ,  framesep=2mm  ,   fontsize=\small]{python}
def is_feasible(self, params):
	x0 = params['x0']
	x1 = params['x1']
	
	if np.logical_or(0.41 < x0 < 0.46, 0.54 < x0 < 0.59):
		return False
	if np.logical_or(0.34 < x1 < 0.41, 0.59 < x1 < 0.66):
	    return False
	return True
\end{minted}


 \item  \textit{Slope} This surface generalizes a plane to discrete domains. The surface's values linearly increase along each dimension. Constraints form three area elements defined by circles with increasing radii.
 
 \begin{minted}[mathescape  ,  autogobble  ,  numbersep=5pt  ,   frame=lines  ,  framesep=2mm  ,   fontsize=\small]{python}
def is_feasible(params):
    x0 = params['x0'] 
    x1 = params['x1']  
    y = x0**2 + x1**2
    if 5 < y < 25:
        return False
    if 70 < y < 110:
        return False
    if 200 < y < 300:
        return False
    return True
\end{minted}

 \item \textit{Sphere}: This surfaces generalizes a parabola to discrete spaces. It features a degenerate global minimum if the number of options along at least one dimension is even,and a well-defined minimum if the number of options for all dimensions is odd. Constraints remove the same two integer inputs, 9 and 11, from consideration in both dimensions.
 
\begin{minted}[mathescape  ,  autogobble  ,  numbersep=5pt  ,   frame=lines  ,  framesep=2mm  ,   fontsize=\small]{python}
def is_feasible(params):
    x0 = params['x0']
    x1 = params['x1']  
    if x0 in [9, 11]:
        return False
    if x1 in [9, 11]:
        return False
    return True
\end{minted}

 \item \textit{Michalewicz}: This surface features a sharper well where the global optimum is located. The number of psuedo-local minima scales factorially with the number of dimensions. Constraints consist of the area element between a circle centred around $(x_0, x_1) = (14, 10)$ with radii $\sqrt{5}$ and $\sqrt{30}$, as  well as two rectangular areas.
 
\begin{minted}[mathescape  ,  autogobble  ,  numbersep=5pt  ,   frame=lines  ,  framesep=2mm  ,   fontsize=\small]{python}
def is_feasible(params):
    x0 = params['x0'] 
    x1 = params['x1'] 
    y = ((x0-14))**2 + (x1-10)**2
    if 5 < y < 30:
        return False
    if 12.5 < x0 < 15.5:
        if x1 < 5.5:
            return False
    if 8.5 < x1 < 11.5:
        if x0 < 9.5:
            return False
    return True
\end{minted}

 \item  \textit{Camel}: This surface features a degenerate and pseudo-disconnected global minimum. In 2-dimensions, it has global minima at $(x_0, x_1) = (7, 11)$ and $(x_0, x_1) = (14, 10)$. Constraints are generated by randomly sampling 100 infeasible locations and excluding the $(x_0, x_1) = (7, 11)$ optima.
 
 \begin{minted}[mathescape,autogobble,numbersep=5pt,frame=lines,framesep=2mm,   fontsize=\small]{python}
def is_feasible(params):
    # choose infeasible points at random
    num_opts = 21
    options = [i for i in range(0,num_opts,1)]
    num_infeas = 100
    np.random.seed(42)
    infeas_arrays = np.array([np.random.choice(options, size=num_infeas,replace=True),  
                              np.random.choice(options, size=num_infeas,   replace=True)]).T
    infeas_tuples = [tuple(x) for x in infeas_arrays]
    
    # always exclude the other minima
    infeas_tuples.append((7, 11))
    infeas_tuples.append((7, 15))
    infeas_tuples.append((13, 5))

    x0 = params['x0']
    x1 = params['x1']  
    sample_tuple = (x0, x1)
    if sample_tuple in infeas_tuples:
        return False
    return True

\end{minted}

\end{itemize}

\subsection{Results of the constrained optimization benchmarks} \label{sisubsubsec:full_synth_cat_res}

\begin{figure*}[!htb]
    \centering
    \includegraphics[width=1.0\columnwidth]{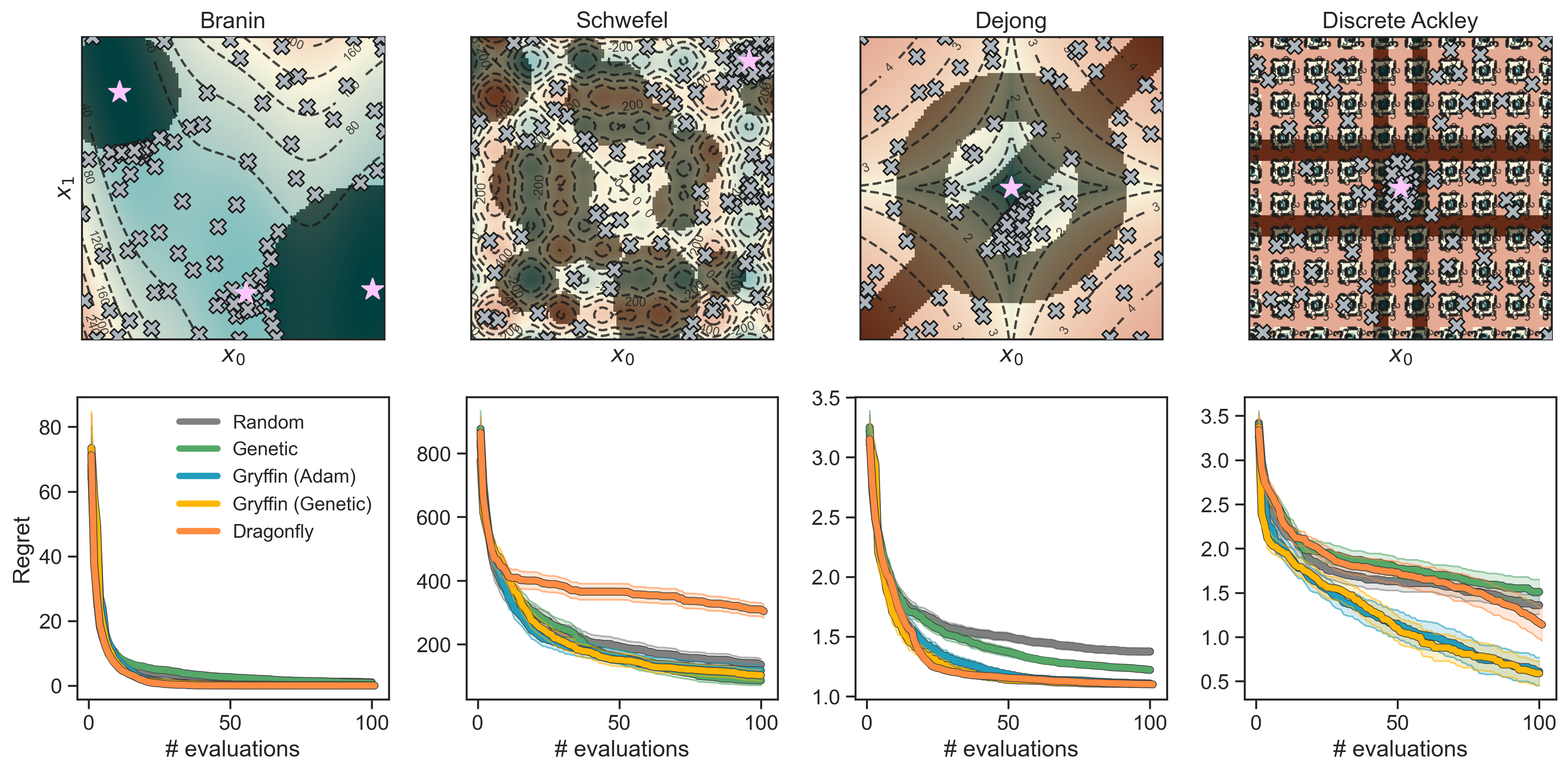}
    \caption{Constrained optimization benchmarks on analytical functions with continuous parameters. The upper row shows contour plots of the surfaces with constrained regions darkly shaded. Gray crosses show sample observation locations and purple stars denote the location(s) of unconstrained global optima. The bottom row show optimization traces for each strategy. Shaded regions around the solid trace represent 95\% confidence intervals.} 
    \label{sifig:continuous_benchmarks_linear}
\end{figure*}

Fig.~\ref{sifig:continuous_benchmarks_linear} shows the results of the continuous optimization benchmarks where regret is displayed on a linear scale, which highlights how performance differences between \gryffin and \dragonfly on \textit{Branin} and \textit{Dejong} are marginal. Table~\ref{sitab:cat_synth_opt_performances} reports the optimization performance achieved by the strategies tested on the discrete surfaces.

\begin{table*}[!htb] 
\begin{center} 
\begin{ruledtabular}
\begin{tabular}{ccccc}
 $-$ & \textbf{Slope} (311) & \textbf{Sphere} (362) & \textbf{Michalewicz} (323)  & \textbf{Camel} (347) \\
\hline
Random  & $157.3 \pm 9.4$ &  $162.3 \pm 9.7$  & $167.8 \pm 9.2$ & $171.0 \pm 10.8$ \\
Genetic & $55.2 \pm 2.7$ &  $61.5 \pm 2.7$  &  $47.7 \pm 2.4$  &  $92.9 \pm 6.2$ \\
\gryffin (Hill)  & $12.7 \pm 1.0$ & $19.0 \pm 0.8$  & $18.4 \pm 0.9$ &  $33.8 \pm 1.6$ \\
\gryffin (Genetic)  & $12.4 \pm 1.0$ & $20.2 \pm 0.8$  & $18.7 \pm 0.8$ & $36.0 \pm 2.6$ \\
\dragonfly & $11.0 \pm 0.1$ &  $13.6 \pm 0.3$  & $29.8 \pm 1.2$ &  $39.0 \pm 2.3$ \\
\end{tabular}
\end{ruledtabular}
\caption{Mean and standard error of the number of evaluations needed by each strategy to identify the global optimum of each constrained discrete surface tested. The integer in parentheses in the header is the number of feasible tiles for the surface after the constraint is applied, out of a total of $21\times21=442$ input combinations.}
\label{sitab:cat_synth_opt_performances}
\end{center}
\end{table*}

\subsection{Empirical comparison of sampling in Gryffin and Dragonfly} \label{sisubsubsec:gryffin_dragonfly_comparison}

\begin{figure*}[!htb]
    \centering
    \includegraphics[width=1.0\columnwidth]{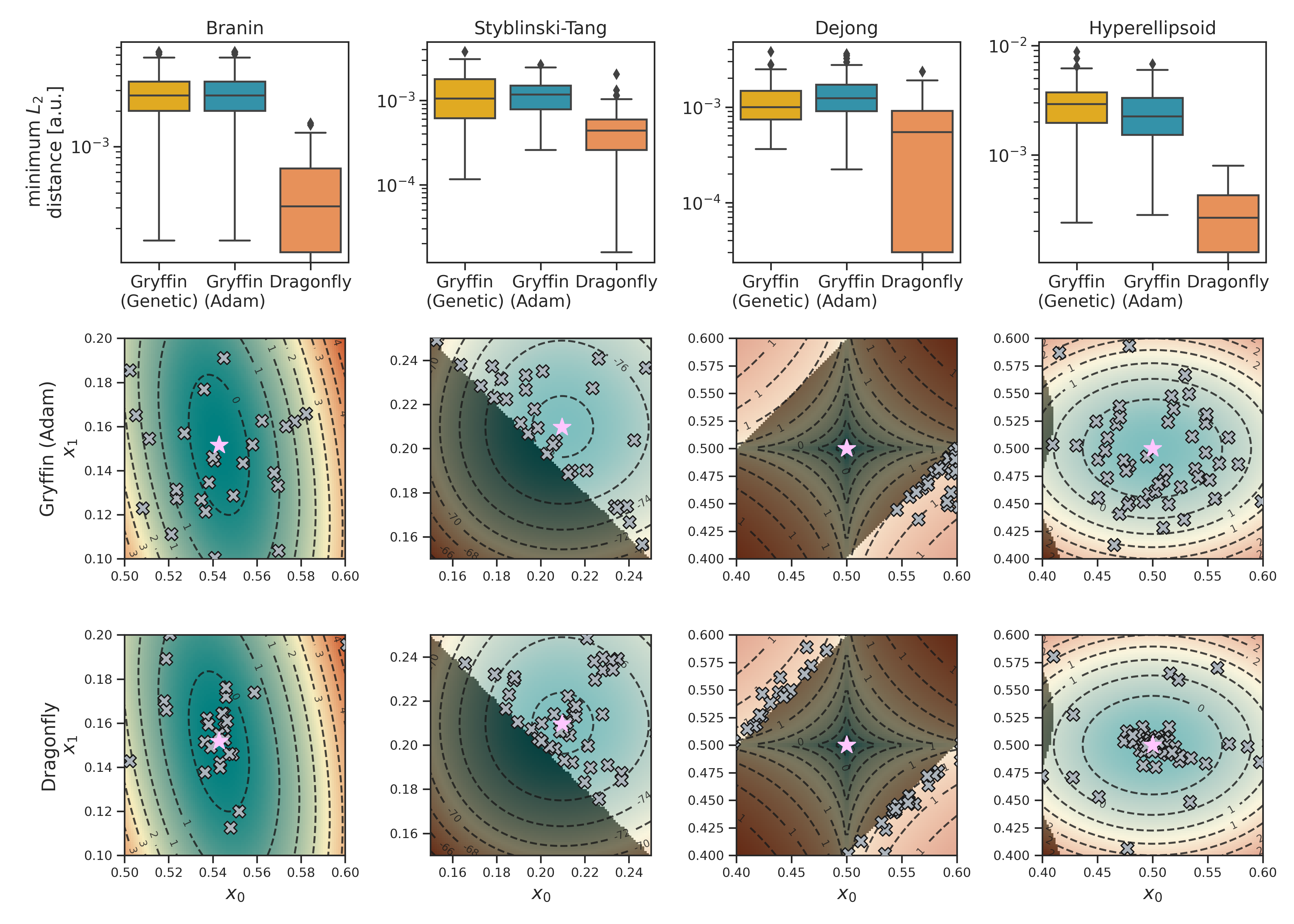}
    \caption{Empirical evaluation of the sampling behaviour of \textit{Gryffin} and \textit{Dragonfly} on constrained continuous surfaces. The first row shows the minimum Euclidean distance between any two parameter points selected by the each optimization strategy. For each continuous constrained surface, \textit{Dragonfly} allows for recommendation of parameter points which are significantly closer to past observations than does \textit{Gryffin (Adam)} or \textit{Gryffin (Genetic)}. The second and third rows shows the location of \textit{Gryffin (Adam)} and \textit{Dragonfly} samples (grey crosses) in the vicinity of the surface minima (pink star).}  
    \label{sifig:gryffin_dragonfly_comparsion_boxplots}
\end{figure*}

In this section, we examine the sampling tendencies of \gryffin and \dragonfly on the constrained, continuous analytical benchmark functions. Specifically, we compare the tendency of each algorithm to suggest parameter point which are in close proximity to past observations. 
The first row of Fig.~\ref{sifig:gryffin_dragonfly_comparsion_boxplots} shows the minimum Euclidean distance between any two parameter points selected during an optimization campaign by each planner (boxplots show this metric over the 100 independently seeded runs). \textit{Dragonfly} is able to recommend parameter points which are significantly closer to past observations than \textit{Gryffin (Adam)} or \textit{Gryffin (Genetic)}. The greater exploitative tendency of \textit{Dragonfly} is beneficial on smooth continuous surfaces as it allows for marginal improvement on regret values (main text Fig.~\ref{fig:synth_continuous}). \textit{Gryffin} strategies, on the other hand, contain a self-avoidance routine which biases the search away from past observations in an attempt to avoid redundant measurements. For practical experimental applications in chemistry, the resolution on input parameters is determined by precision of laboratory equipment and/or human error, and should be considered before commencing the experiment. 
The bottom two rows show the location of observations for \textit{Gyrffin (Adam)} and \textit{Dragonfly} strategies around the minima of each surface. Visually, it is apparent that \textit{Dragonfly} has a greater tendency to recommend parameter points which are considerably closer to past observations than does \textit{Gryffin}.

\section{Process-constrained optimization of \textit{o}-xylenyl C\textsubscript{60} adducts synthesis} \label{sisubsec:fullerene_details}

\subsection{Details of the Bayesian neural network experiment emulator} \label{sisubsubsec:fullerene_emulator}

\begin{figure*}[!ht]
    \centering
    \includegraphics[width=1.0\columnwidth]{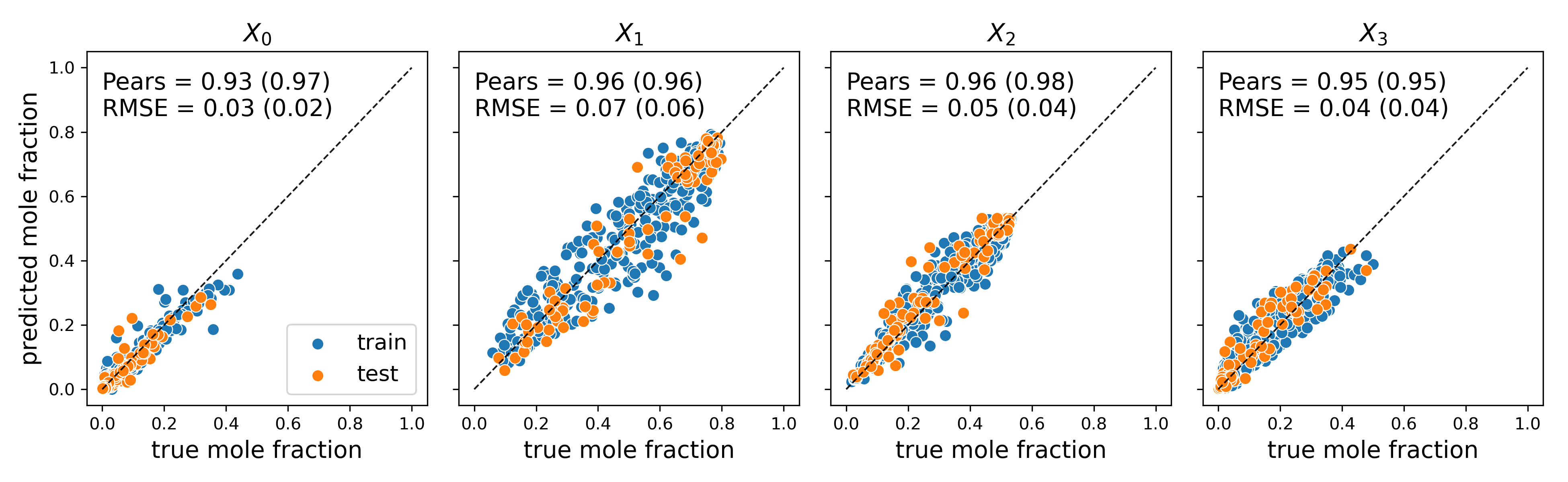}
    \caption{ Parity plots for each mole fraction predicted by our Bayesian neural network emulator, averaged over 50 network parameter samples. Horizontal axes plot the true mole fraction, and vertical axes plot the predicted mole fraction. The dashed diagonal line indicates perfect agreement. The Pearson correlation coefficient and root-mean-square error is given for the training set (in parentheses) and test set for each mole fraction target. Train (test) set points are shown in blue (orange). } 
    \label{sifig:emulator_performance}
\end{figure*}

To emulate the process-constrained synthesis of C\textsubscript{60} adducts, we trained a Bayesian neural network (BNN) to return stochastic outcomes based on a set of controllable parameters. The trained emulator takes a vector containing the experimental conditions ($T$, $F_{\text{C}_{60}}$, and $F_\text{S}$) and predicts the mole fractions of the products, the un- ($[X_0]$), singly- ($[X_1]$), doubly- ($[X_2]$), and triply-functionalized ($[X_3]$) $\text{C}_{60}$. The BNN consisted of 3 densely-connected variational layers with reparameterized Monte Carlo estimators\cite{blundell2015weight} and was implemented in PyTorch.\cite{NEURIPS2019_9015}. Each hidden layer had $64$ nodes and featured a ReLU non-linearity, while the output layer had $4$ nodes, one for each of the aformentioned $\text{C}_{60}$ adducts. The output layer used the softmax activation function, which normalizes the outputs to a probability distribution where $\sum_{i=0}^3 [X_i] = 1$. Network weights $w_i$ and biases $b_i$ followed Gaussian distributions whose priors were set to have zero mean and unit standard deviation, i.e. $w_i \sim \mathcal{N}(0,1)$,   $b_i \sim \mathcal{N}(0,1)$. The network was trained using variational Bayesian inference. The ELBO loss was minimized using the Adam optimizer\cite{kingma_adam_2017}, and resulting gradients were used to adjust the weights and biases of the network's parameter distributions during training. Fig.~\ref{sifig:emulator_performance} shows parity plots of our model's predictions against the true $\text{C}_{60}$ adduct mole fractions, where $500$ experimental measurements were used for training and $100$ for testing. The BNN displayed excellent interpolation performance across the parameter space for each adduct type, with Pearson correlation coefficients on the test sets between $0.93$ and $0.96$.

\subsection{Estimating the experimental cost} \label{sisubsubsec:fullerene_cost_estimation}

The overall goal of the process-constrained optimization of \textit{o}-xylenyl C\textsubscript{60} adducts synthesis is to adjust reaction conditions such that the combined yield of first- and second-order adducts is maximized and reaches at least 90\%, while the cost of reagents is minimized. In order to estimate the cost of the experiments, we considered the listed price of sultine and C\textsubscript{60} by the chemical supplier Sigma-Aldrich. The cost of dibromo-o-xylene cost on Sigma Aldrich was $\$191$ for 100 g. The cost of C\textsubscript{60} was $\$422$ for 5g. In the experiments by Walker \textit{et al.}~\cite{walker_tuning_2017}, the concentration of sultine was $1.4$ mg/mL, while the concentration of C\textsubscript{60} was $2.0$ mg/mL. The amount of C\textsubscript{60} used in the experiments will therefore have much greater influence on overall experiment cost than sultine. Our optimization experiments target the adjustment of volume flow rates of each of these chemicals. Thus, we seek a measure of per-unit-time operation cost to be minimized. Converting to per-litre costs, we have $2.674 \; \$/\text{L}$ for sultine, and $168.8 \; \$/\text{L}$ for C\textsubscript{60}. Finally, from the flow rates used in the experiment, $F_{\text{C}}$ and $F_{\text{S}}$ (with units of $\mu\text{L}/\text{min}$), we obtain an estimate of per-minute operation cost of the flow-reactor from Walker \textit{et al.}~\cite{walker_tuning_2017} with units of $\$/\text{min}$ as

\begin{align} \label{sieq:cost_expression}
    \text{cost} = \frac{1 \text{L}}{10^6 \mu \text{L}} F_{\text{C}} \times \frac{\$ 168.8}{\text{L}} + \frac{1 \text{L}}{10^6 \mu \text{L}} F_{\text{S}} \times \frac{\$ 2.674}{\text{L}} \,.
\end{align}

\begin{figure*}[htb]
    \centering
    \includegraphics[width=1.0\columnwidth]{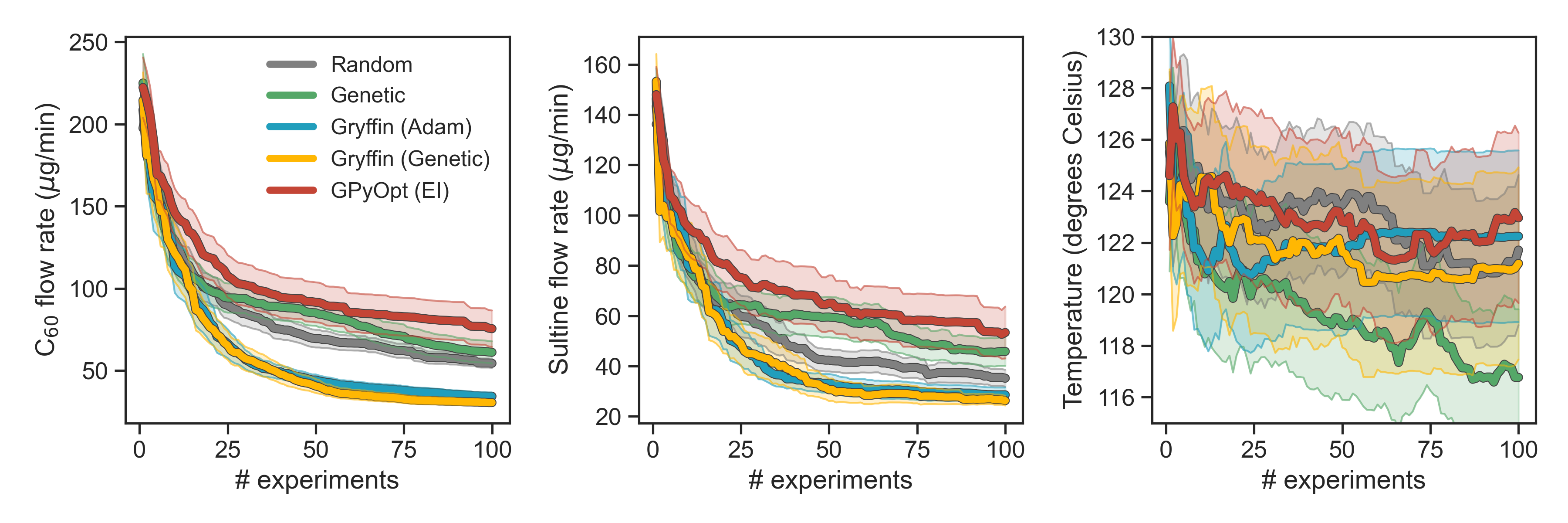}
    \caption{Mean and 95\% confidence interval for parameter values corresponding to the best objective values found by each optimization strategy at each iteration of the optimization campaign. C\textsubscript{60} and sultine flow rates are both shown with units of $\mu$g/min.} 
    \label{sifig:fullerene_100_parameter_traces}
\end{figure*}

Fig.~\ref{sifig:fullerene_100_parameter_traces} shows the mean and 95\% confidence interval for parameter values corresponding to the best observed objective values achieved by each optimization strategy. We report the flow rate in terms of \textit{mass} per unit time (mass flow rate) to account for the difference in concentrations of each reagent and compare the rates on an equal footing. To improve upon our secondary cost objective, each strategy decreases the $F_{C}$ parameter, as it's value dominates the cost in Eq.~\ref{sieq:cost_expression}. Decrease in $F_{C}$ is however accompanied by a decrease in $F_S$ to preserve the high ($\geq 0.9$) mol fractions of the $X_1$ and $X_2$ adducts. For most optimization runs, the temperature of the best performing reactions varies between 116 and 130 $^{\circ}$C.

\begin{figure*}[htb]
    \centering
    \includegraphics[width=0.8\columnwidth]{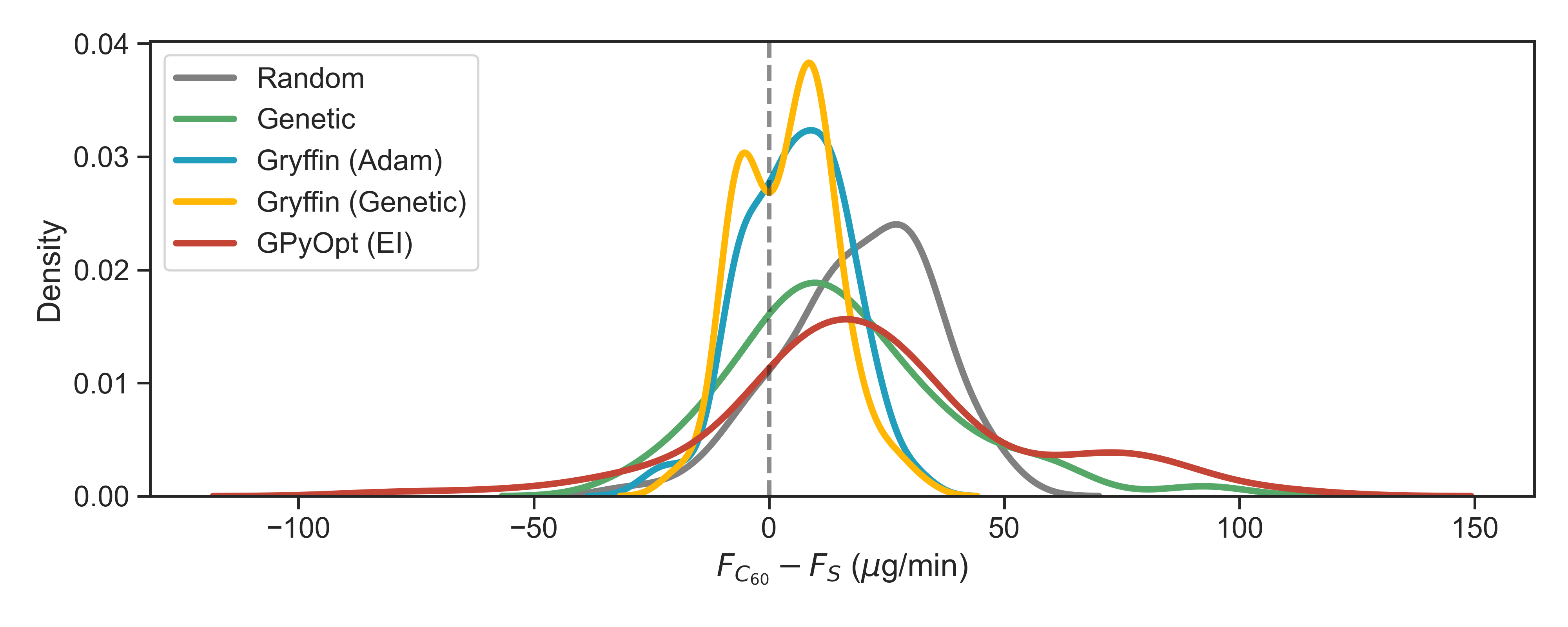}
    \caption{ Kernel density estimates show the distribution of $F_{\text{C}} - F_{\text{S}}$ for the best reactions conditions achieved by each optimization strategy in units of $\mu$g/min. Each distribution is comprised of 100 such values, one for each independently seeded optimization. Positive values indicate that the best achieved reaction conditions had $F_{\text{C}} > F_{\text{S}}$.  } 
    \label{sifig:fullerene_flow_rate_diff_kde}
\end{figure*}

Fig.~\ref{sifig:fullerene_flow_rate_diff_kde} shows distributions of $F_{\text{C}} - F_{\text{S}}$ values for the best reaction conditions achieved by each optimization strategy in units of $\mu$g/min. For each strategy, we note that this distribution favours positive values, meaning that, in the majority of the best achieved reaction conditions, the mass flow rate of  C\textsubscript{60} was greater than that of sultine. Crucially, the \textit{Gryffin} strategies, which exhibited the best optimization performance on this application, achieved narrower distributions around $F_{\text{C}} - F_{\text{S}} = 0$ than other strategies.

\section{Design of redox-active materials for flow batteries with synthetic accessibility constraints} \label{sisubsec:redox_details}

\subsection{Computation of reduction potential tolerance} \label{sisubsubsec:dft_reduction_constraint}

To set the reduction potential ($E^{\text{red}}$) upon which we would like to improve,   and which is used as an absolute tolerance in \chimera\cite{hase_chimera_2018}, we computed $E^{\text{red}}$ for the base scaffold molecule \textbf{H-AcBzC\textsubscript{6}}.\cite{agarwal_discovery_2021} We computed $E^{\text{red}}$ with the same computational protocol used by Agarwal \textit{et al.}~\cite{agarwal_discovery_2021}. The DFT calculation was performed using Gaussian 16~\cite{g16} at  the wb97xd/6-31+G-(d,p)~\cite{chai_long-range_2008,rassolov_6-31g_2001} level of theory. Optimized neutral and anionic geometries were subject to frequency calculations to compute the free energies. The SMD continuum model~\cite{marenich_universal_2009} was used with acetonitrile as the solvent. The reduction potential was calculated using Eq. 1 in Agarwal \textit{et al.}~\cite{agarwal_discovery_2021},  

\begin{align}
    E^{\text{red}} = \frac{-\Delta G^{\text{red}}}{nF} - 1.24 \; \text{V} \,,  
\end{align}

where $\Delta G^{\text{red}} = G^{\text{reduced}} - G^{\text{neutral}} $, $n$ is the number of electrons added to the neutral molecules ($n=1$), $F$ is Faraday's constant in eV, and 1.24 is a constant subtracted to convert the Gibbs free energy change to reduction potential (with Li/Li\textsuperscript{+} reference electrode). The $E^{\text{red}}$ for \textbf{H-AcBzC\textsubscript{6}} was computed to be $2.038372$ V. Of the $1408$ functional derivatives subject to computation by Agarwal \textit{et al.}~\cite{agarwal_discovery_2021}, only $243$ had better (lower) $E^{\text{red}}$.

\subsection{Prediction of the synthetic accessibility of redoxmer candidate molecules} \label{sisubsubsec:battery_synthetic_access}

As a constraint on the redoxmer candidates space, we enforce a retrosynthetic accessibility threshold below which the candidate is considered infeasible. The goal was to have an indication of synthetic accessibility that could be used to constrain the search space to candidates that likely to be synthesizable in practice. 

\begin{figure*}[!ht]
    \centering
    \includegraphics[width=1.0\columnwidth]{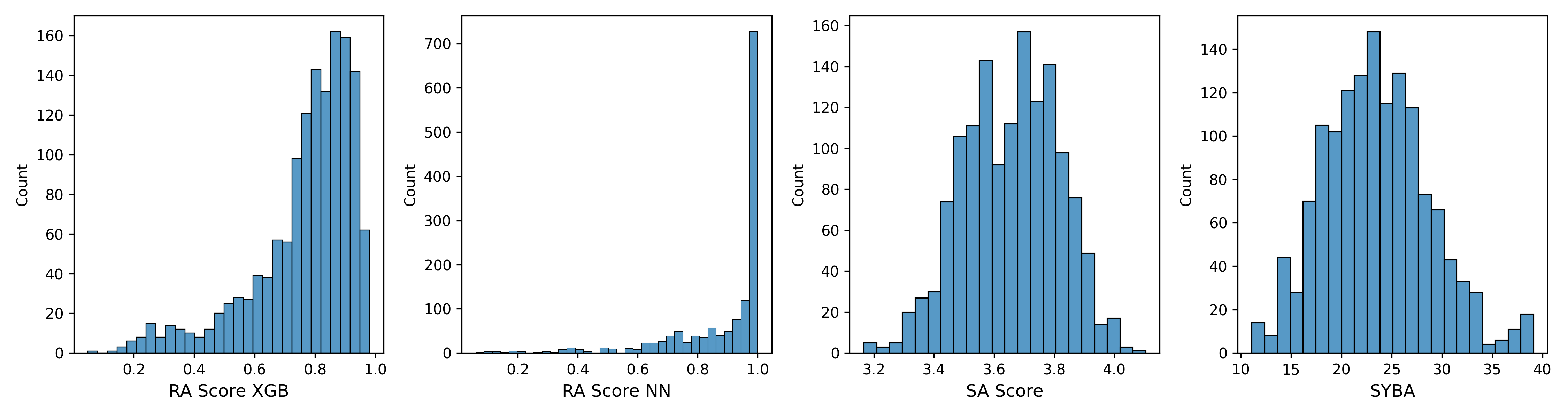}
    \caption{Histograms showing the distributions of four synthetic accessibility scores computed for the $1408$ redoxmer candidates.} 
    \label{sifig:synthetic_access_scores}
\end{figure*}

Fig.~\ref{sifig:synthetic_access_scores} shows the distributions of different synthetic accessibility scores for the set of $1408$ redoxmer candidates considered in this application. Specifically, it includes the \textit{RAscore}~\cite{thakkar_retrosynthetic_2021} predicted by an XGBoost classifier (XGB) and a neural network (NN), the fragment-based synthetic accessibility score \textit{SAscore}\cite{ertl_estimation_2009}, and the synthetic Bayesian classifier (\textit{SYBA})\cite{Vorsilak:2021}. The \textit{RAScore} is a recently reported synthetic accessibility score that tries to capture the probability of \textit{AiZynthFinder} being able to identify a synthetic route for the molecule being evaluated.~\cite{thakkar_retrosynthetic_2021} \textit{AiZynthFinder} is a retrosynthetic planning tool that can generate synthetic routes for organic molecules.\cite{genheden_aizynthfinder_2020} Hence, an \textit{RAScore} of $1$ indicates a synthetic path to the desired molecule is likely to exist, while a score of $0$ indicates that finding a synthetic path is likely to be challenging and potentially impossible. For the purpose of our constrained optimization experiments, we decided to use the \textit{RAscore} based on a NN model given its reported performance\cite{thakkar_retrosynthetic_2021} and intuitive interpretation.

\subsection{Generation of descriptors for benzothiadiazole scaffold substituents} \label{sisubsubsec:battery_descriptors}

In this example application of constrained Bayesian optimization, we employed the \textit{Dynamic} version of \gryffin for combinatorial optimization,\cite{hase_gryffin_2021} which can take advantage of physicochemical descriptors in the search for optimal molecules. Specifically, we provided \gryffin with a total of seven simple descriptors associated with each of the four substituent groups considered ($R_{1-4}$ in Fig.~\ref{fig:batteries_multipanel}a). The physicochemical descriptors were computed with the \textit{Mordred} Python package.~\cite{moriwaki_mordred_2018} As summarized in Table~\ref{sitab:r_group_descriptors_objectives}, the following descriptors were considered: the number of hetero atoms (\texttt{nHetero}), molecular weight (\texttt{MW}), topological polar surface area (\texttt{TopoPSA}), number of heavy atoms (\texttt{nHeavyAtom}), atomic polarizablity (\texttt{apol}), fraction of sp\textsuperscript{3} hybridized carbons (\texttt{FCSP3}), and geometric diameter (\texttt{Diameter}). All seven descriptors were used for substituent groups $R_{2-4}$, but only four of them are used for $R_1$. We eliminate \texttt{apol}, \texttt{FCSP3} and \texttt{Diameter} from consideration because they each have equal value for both $R_1$ substituent options, and therefore are not informative. Table~\ref{sitab:r_group_descriptors_objectives} sumarizes the Pearson correlation of each descriptor with each objective value over the entire set of 1408 molecules.  N/A entries show the cases where the descriptor is omitted for the $R_1$ substituent. In addition, Table~\ref{sitab:r_group_descriptors_objectives} reports the Pearson correlations between the descriptors for all four R-groups ($\rho_1$, $\rho_2$, $\rho_3$, and $\rho_4$) and each optimization objective ($\Delta \lambda^{\text{abs}}$, $E^{\text{red}}$, and $G^{\text{solv}}$). Table~\ref{sitab:r_group_descriptors_pairwise} reports instead the pairwise correlation between each descriptor, averaged over all R-groups.

\begin{table*}[htb] 
 \begin{center} 
\begin{ruledtabular}
\begin{tabular}{c c c c c c c c c c c c c}
Mordred name & \multicolumn{4}{c}{$\Delta \lambda^{\text{abs}}$}  & \multicolumn{4}{c}{$E^{\text{red}}$} &  \multicolumn{4}{c}{$G^{\text{solv}}$}  \\ 
 & $\rho_1$ & $\rho_2$ & $\rho_3$ & $\rho_4$ & $\rho_1$ & $\rho_2$ & $\rho_3$ & $\rho_4$ & $\rho_1$ & $\rho_2$ & $\rho_3$ & $\rho_4$   \\ 
\hline
\texttt{nHetero}   & \textbf{0.17}  &  0.13  &   0.12  &   0.19    &  \textbf{0.22}  &    0.27  &   0.26  &   0.35                    & \textbf{0.62}  &  \textbf{0.24}  &   \textbf{0.22}  &  -0.12 \\
\texttt{MW}        &  \textbf{0.17}  &   0.06  &   0.05  &   0.15                   &  \textbf{0.22}  &   0.21  &   0.20  &   0.19                                         & \textbf{0.62}  &   0.23  &   0.21  &   -0.13 \\ 
\texttt{TopoPSA}   & \textbf{-0.17}  &   -0.12  &   -0.11  &   -0.14               & \textbf{-0.22}  &   0.13  &   0.11  &   0.30                                                  &  \textbf{-0.62}  &   0.00  &   0.01  &   -0.22  \\
\texttt{nHeavyAtom} &  \textbf{0.17}  &  0.09  &  0.07  &  \textbf{0.20}             & \textbf{0.22}  &   0.22  &   0.21  &   0.22                                          &  \textbf{0.62}  &   0.22  &   0.20  &   -0.21  \\
\texttt{apol}  & \textbf{-0.17}  &   \textbf{-0.25}  &   \textbf{-0.31}  &   -0.14     & \textbf{-0.22}  &   -0.26  &  -0.27  &   -0.16                                                &  \textbf{-0.62}  &   -0.06  &   -0.04  &   \textbf{-0.41}  \\
\texttt{FCSP3}  & N/A  &   -0.11  &   -0.06  &   -0.04                           &  N/A  &   \textbf{-0.31}  &   \textbf{-0.32}  &   \textbf{-0.36}                                                  &  N/A  &   0.13  &   0.14  &   0.02 \\
\texttt{Diameter} & N/A  &  0.02  &  -0.04  &  0.17                                & N/A  &   0.16  &   0.16  &   0.23                                                     & N/A  &   0.14  &   0.12  &   -0.29 \\
\end{tabular}
\end{ruledtabular}
\caption{Mordred descriptors used to describe R-groups for the battery application optimization. The right most three columns show the Pearson correlation between each descriptor and each optimization objective for the $1408$ redoxmer candidates considered. The correlations for each of the four R groups are comma separated  ,   i.e. $\rho_{R1}  ,  \rho_{R3}  ,  \rho_{R4}  ,  \rho_{R5}$. The largest correlation for each objective and R-group is bolded. N/A entries indicate that this descriptor was not considered for this particular $R$ group. For the $R_1$ group, we do dont consider \texttt{apol}, \texttt{FCSP3} and \texttt{Diameter} since their values are the same for both $R_1$ options and therefore provide no additional information.}
\label{sitab:r_group_descriptors_objectives}
\end{center}
\end{table*}

\begin{table*}[htb] 
 \begin{center} 
\begin{ruledtabular}
\begin{tabular}{cccccccc}
 $-$ & \texttt{nHetero} & \texttt{MW} & \texttt{TopoPSA}  & \texttt{nHeavyAtom} & \texttt{apol} & \texttt{FCSP3} &  \texttt{Diameter}  \\ 
\hline
\texttt{nHetero}  & $1.00$ & $0.92$  & $0.01$ & $0.91$ & $-0.04$  & $0.2$ & $0.61$   \\
\texttt{MW}    & $0.92$ & $1.00$  & $-0.06$  & $0.93$  & $0.24$  & $0.29$  & $0.71$ \\
\texttt{TopoPSA}    & $0.01$  & $-0.06$  & $1.00$ & $-0.14$   & $0.03$  & $-0.20$  &  $-0.13$ \\
\texttt{nHeavyAtom}    & $0.91$  & $0.93$  & $-0.14$   &$1.00$  & $0.32$  & $0.35$ & $0.86$\\
\texttt{apol}    & $-0.04$  & $0.24$   & $0.03$   &  $0.32$ & $1.00$ & $0.47$ & $0.62$ \\
\texttt{FCSP3}    & $0.2$ & $0.29$  & $-0.2$  & $0.35$  & $0.47$  & $1.00$ & $0.32$ \\
\texttt{Diameter}    & $0.61$  & $0.71$   & $-0.13$   & $0.86$  & $0.62$  & $0.32$  & $1.00$ \\
\end{tabular}
\end{ruledtabular}
\caption{Pairwise Pearson correlations between Mordred descriptors used to describe the R-groups of the redoxmer candidates.}
\label{sitab:r_group_descriptors_pairwise}
\end{center}
\end{table*}

\subsection{Additional optimization experiments}
\begin{figure*}[htb]
    \centering
    \includegraphics[width=1.0\columnwidth]{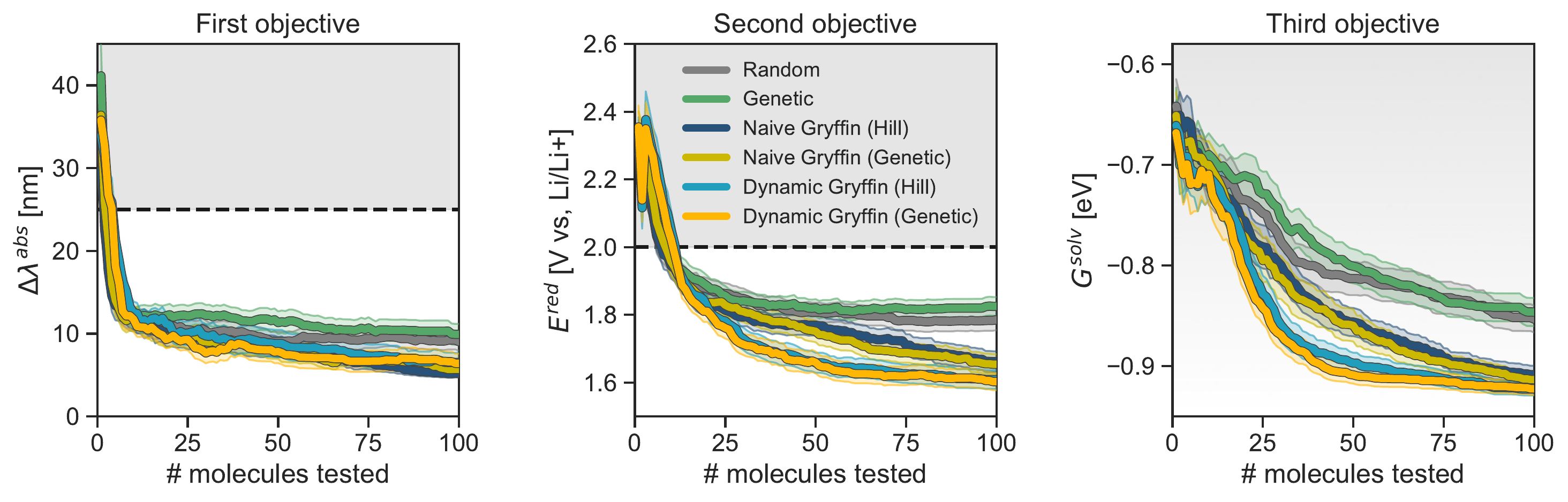}
    \caption{Results of the constrained optimization experiments for the design of redox-active flow battery materials. Grey shaded regions indicate objective values failing to achieve the desired objectives. Traces depict the objective values corresponding to the best achieved merit at each iteration, where error bars represent 95\% confidence intervals.} 
    \label{sifig:batteries_traces}
\end{figure*}

In addition to \gryffin optimizations taking advantage of physicochemical descriptors (\textit{Dynamic Gryffin}), we also carried out optimizations without this additional information using \textit{Naive Gryffin}. Fig.~\ref{sifig:batteries_traces} shows the optimization performance of all strategies tested, including the latter. The results show how the use of descriptors provide an edge to \gryffin to achieve superior performance to all other strategies. Regardless, \textit{Naive Gryffin} still outperforms model-free optimization strategies \textit{Random} and \textit{Genetic}. All these optimizations were constrained to molecules with high synthetic accessibility scores, as described above.

	\putbib[main]
\end{bibunit}




\end{document}